\documentclass[11pt]{article}
\usepackage{amsthm,amssymb,amscd,amsmath,latexsym,indentfirst,color,amsfonts}
\usepackage[mathcal]{eucal}

\theoremstyle{plain}
\newtheorem{theorem}{Theorem}[section]
\newtheorem{prop}[theorem]{Proposition}
\newtheorem{coro}[theorem]{Corollary}
\newtheorem{lem}[theorem]{Lemma}
\newtheorem{defi}[theorem]{Definition}

\theoremstyle{remark}
\newtheorem{rem}[theorem]{Remark}
\newtheorem{ex}[theorem]{Example}

\DeclareMathOperator\PGL{PGL}
 \DeclareMathOperator\Aut{Aut}
\DeclareMathOperator\Hom{Hom}\DeclareMathOperator\Int{Int}

\DeclareMathOperator\an{and}
\DeclareMathOperator\Stab{Stab}\DeclareMathOperator\Rad{Rad}
\DeclareMathOperator\Supp{Supp}

\def\a{\frak a}
\def\sl{\frak s\frak l} \def\sp{\frak s\frak p}  \def\so{\frak s\frak o}
\def\gl{\frak g\frak l}

\def\bt{\beta}
\def\si{\sigma}\def\ve{\varepsilon}
\def\al{\alpha}\def\Ga{\Gamma}
\def\w{\wedge}\def\rt{\rightarrow}\def\st{\subset}\def\ot{\otimes}\def\op{\oplus}
\def\wh{\widehat}

\newcommand{\be}{\begin {equation}}
\newcommand{\ee}{\end {equation}}
\newcommand{\bp}{\begin {proof}}
\newcommand{\ep}{\end {proof}}
\newcommand{\bee}{\begin {equation*}}
\newcommand{\eee}{\end {equation*}}
\newcommand{\lb}{\label}

\DeclareMathOperator\U{U}

\def\O{\operatorname{O}}

\DeclareMathOperator\PO{PO}

\DeclareMathOperator\PSp{PSp}

\def\diag{\operatorname{diag}}
\def\E{\operatorname{E}}
\newcommand{\bbR}{\mathbb{R}}
\newcommand{\fre}{\mathfrak{e}}
\newcommand{\frf}{\mathfrak{f}}
\newcommand{\frg}{\mathfrak{g}}

\def\a{\frak a}

\def\sl{\frak s\frak l}
\def\sp{\frak s\frak p}

\def\so{\frak s\frak o}
\def\gl{\frak g\frak l}

\begin{document}

\title{Fine gradings of complex simple Lie algebras and Finite Root Systems}

\author{Gang Han \\ Department of Mathematics \\ Zhejiang University\\China \thanks{Corresponding author. Work is supported by Zhejiang Province Science Foundation, grant No. LY14A010018.}
\and Kang Lu \\ Department of Mathematics \\ Zhejiang University\\China
\and Jun Yu \\ Department of Mathematics \\ Massachusetts Institute of Technology \\ Cambridge, MA 02139, USA }

\date{October 29, 2014}


\maketitle

{\small \noindent \textbf{Abstract}.

A $G$-grading on a complex semisimple Lie algebra $L$, where $G$ is a finite abelian group, is called quasi-good if
each homogeneous component is 1-dimensional and 0 is not in the support of the grading.

Analogous to classical root systems, we define a finite root system $R$ to be some subset of a finite symplectic abelian group satisfying certain axioms. There always
corresponds to $R$ a semisimple Lie algebra $L(R)$ together with a quasi-good grading on it. Thus one can construct nice basis of $L(R)$ by means of finite root systems.

We classify finite maximal abelian subgroups $T$ in $\Aut(L)$ for complex simple Lie algebras $L$ such that the grading induced by the action of $T$ on $L$ is quasi-good, and show that the set of roots of $T$ in $L$ is always a finite root system. There are five series of such finite maximal abelian subgroups, which occur only if $L$ is a classical simple Lie algebra.

}

\tableofcontents

\section{Introduction}
 \setcounter{equation}{0}\setcounter{theorem}{0}

\noindent \textbf{1.1.} We will first briefly review gradings on a (not necessarily associative) algebra by abelian groups and mainly focus on gradings on a Lie algebra. All the algebras in the paper are assumed to be finite dimensional and over the field $\mathbb{C}$ of complex numbers, although all of our main results can be generalized to any algebraically closed field of characteristic 0.

Let $A$ be a finite dimensional algebra and $G$ be a finitely generated additive abelian group. A $G$-\textit{grading} $\Gamma$ on $A$ is the decomposition of
$A$ into a direct sum of subspaces $$\Gamma: A=\oplus_{g\in G}
\ A_g$$ such that
$$A_g\cdot A_h\st A_{g+h},\ \forall\ g,h\in G.$$ We also say that $A$ is $G$-\textit{graded}. If $A$ is a Lie algebra then the multiplication is understood to be the Lie bracket. The subset
$R=\{g\in G|A_g\neq 0\}$ of $G$ is called the \textit{support} of this grading and is denoted $\Supp\ \Ga$. For any $g\in R$, $A_g$ is called the \textit{homogeneous component} of degree $g$ and each element in $A_g$ is a \textit{homogeneous} element. We will always
assume that $G$ is generated by $R$, otherwise it could be
replaced by its subgroup generated by $R$. So $G$ is always a finitely generated abelian group. Denote this grading by $(\Ga,G)$, or simply by $\Ga$. All the gradings mentioned in the paper will be some $G$-grading with $G$ a (finitely generated) additive abelian group.

Let $(\Ga_i,G_i)$ be two gradings on $A$ for $i=1,2$. Assume $\phi\in \Aut(A)$.  If there is a group isomorphism $\psi:G_1\rt G_2$ such that $\psi(\Supp \ \Ga_1)=\Supp \ \Ga_2$ and for any $a\in \Supp\ \Ga_1$, $\phi(A_a)=A_{\psi(a)}$ , then $\phi$ is called a \textit{grading homomorphism} from
$(\Ga_1,G_1)$ to $(\Ga_2,G_2)$ and denoted $\phi:(\Ga_1,G_1)\rt(\Ga_2,G_2)$. The grading homomorphism $\phi:(\Ga_1,G_1)\rt(\Ga_2,G_2)$ is an \textit{isomorphism} of gradings (which is called \textit{group equivalence} in \cite{k}) if $\phi^{-1}$ is also a grading homomorphism from
$(\Ga_2,G_2)$ to $(\Ga_1,G_1)$. In this case the corresponding group homomorphism is an isomorphism.

Now let $A$ be a Lie algebra and denote it by $L$. Let $\Aut(L)$ (resp. $\Int(L)$) be respectively the group of automorphisms (resp. inner automorphisms) of $L$. Let the group $K$ be either $\Aut(L)$ or $\Int(L)$.
Let $\Ga$ be a $G$-grading on $L$. The group $\Aut_K(\Ga)$ consists of all the automorphisms in $K$ that permute the homogeneous components of $\Ga$. The group $\Stab_K(\Ga)$ consists of all the automorphisms in $K$ that stabilize each homogeneous component of $\Ga$.  Clearly $\Stab_K(\Ga)$ is a normal subgroup of $\Aut_K(\Ga)$, and let $$W_K(\Ga)=_{def}\Aut_K(\Ga)/\Stab_K(\Ga)$$ be the \textit{Weyl} group of the grading $\Ga$ with respect to $K$, which describes the symmetry of $\Ga$. Note that the Weyl group of $\Ga$ defined in \cite{ek} is with respect to $K=\Aut(L)$.

People have done a lot of work in the area of abelian group gradings on Lie algebras (and on other types of algebras) and find many important applications for them. A good survey for it can be found in \cite{k}. The classical Cartan decomposition of a semisimple Lie algebra $L$ is in fact a grading on $L$ by a free abelian group of the same rank as that of $L$. Alekseevskii found a class of interesting gradings on simple Lie algebras and computed the corresponding Weyl groups in \cite{a}. Patera and his collaborators began to study abelian group gradings on Lie algebras systematically in \cite{pz} and \cite{hpp}, and they classified most of the fine gradings on classical simple Lie algebras. Later in \cite{bsz} Bahturin and his collaborators described all the abelian group gradings on classical simple Lie algebras except the simple Lie algebra of type $\textbf{D}_4$ (and on a class of simple Jordan algebras as well). In \cite{bk} Bahturin and Kochetov classified the isomorphism classes of abelian group gradings on all the simple Lie algebras except the simple Lie algebra of type $\textbf{D}_4$ in terms of numerical and group-theoretical invariants. Draper and his collaborators classified all the fine gradings on the simple Lie algebra of type $\textbf{D}_4$ in \cite{dv}, classified all the gradings on the simple Lie algebra of type $\textbf{G}_2$ in \cite{dm1} and classified all the fine gradings on the simple Lie algebra of type $\textbf{F}_4$ in \cite{dm2}. In \cite{ek} Elduque and Kochetov   computed the Weyl groups for all the fine gradings on all the classical simple Lie algebras except $\textbf{D}_4$.\bigskip

\noindent \textbf{1.2.} Now we recall the duality between abelian group gradings and abelian group actions.

An abelian algebraic group consisting of semisimple elements is called a quasitorus. Assume $T$ is a complex quasitorus. Let $\wh{T}$ be the set of algebraic group homomorphisms from $T$ to $\mathbb{C}^\times$. Then $\wh{T}$ is an abelian group, called the \textit{character group} or \textit{dual group} of $T$ with addition defined by
\[ (\al+\bt)(a)=\al(a)\cdot\bt(a),\ for\ any\ \al,\bt\in \wh{T}; \ for\ any\ a\in T. \]The group $\Aut(A)$ is a linear algebraic group. Assume $T\st\Aut(A)$ is a quasitorus. Then the action of $T$ on $A$ induces a $\wh{T}$-grading on $A$:
\be \Ga:A=\op_{g\in \wh{T}}A_g. \ee Each element $g$ in the support of this grading is called a \textit{root} of $T$ in $A$, and $A_g$ is called the \textit{root space} of $g$. If $T^{'}$ is another quasitorus in $\Aut(A)$, then $T$ and $T^{'}$ induces isomorphic group gradings on $A$ if and only if $T$ and $T^{'}$ are conjugate in $\Aut(A)$.

Conversely given a grading $A=\oplus_{g\in G}\ A_g$ with $G$ a finitely generated abelian group, let $\wh{G}$ be the dual group of $G$, which is a  (multiplicative) abelian group consisting of homomorphisms from $G$ to $\mathbb{C}^\times$. Then the grading induces a $\wh{G}$-action on $A$:
\be \si\cdot X=\si(g)X,\ for\ all\ \si\in G, X \in A_g.\ee This action gives a homomorphism $\wh{G}\rt \Aut(A)$, which is injective as the support of the grading generates $G$. So the homomorphism embeds $\wh{G}$ as a quasitorus in $\Aut(A)$ and we will identify $\wh{G}$ with its image in $\Aut(A)$.

Assume that $A$ is a semisimple Lie algebra and denote it by $L$. Assume $T\st\Aut(L)$ is a quasitorus. Then the action of $T$ on $L$ induces a $\wh{T}$-grading $\Ga$ on $L$. If  $T$ is a maximal quasitorus in $\Aut(L)$ (resp. in $\Int(L)$) then $\Ga$ is called a \textit{fine grading} (resp. \textit{fine inner grading}). A grading $\Ga$ is \textit{quasi-good} if 0 is not in $\Supp\ \Ga$ and each homogeneous component is 1-dimensional. A quasitorus $T\st \Aut(L)$ induces a quasi-good grading on $L$ if and only if $T$ is a finite quasitorus in $\Aut(L)$ such that its centralizer in $\Aut(L)$ is also finite. If  $T$ is a finite maximal quasitorus in $\Aut(L)$ (resp. in $\Int(L)$) such that the induced grading $\Ga$ on $L$ is quasi-good, then $\Ga$ is called a \textit{good grading} (resp. \textit{good inner grading}) and $T$ is called a \textit{good} finite maximal quasitorus in $\Aut(L)$ (resp. in $\Int(L)$).

 A good grading must be quasi-good, but not conversely. See Lemma \ref{x1} for details.

Let $L$ be a simple Lie algebra and the group $K$ be either $\Aut(L)$ or $\Int(L)$. Let $T\st K$ be a maximal quasitorus, then it induces a $G$-grading $\Ga$ on $L$ with $G=\wh{T}$.  Define the Weyl group of $T$ with respect to $K$ to be $$W_K(T)=N_K(T)/Z_K(T).$$  By Proposition 2.4 of \cite{h1}, one has $N_K(T)=\Aut_K(\Ga)$, $Z_K(T)=\Stab_K(\Ga)$ and thus \be\lb{l1}W_K(T)=W_K(\Ga).\ee\bigskip

\noindent \textbf{1.3.} Given a finite maximal abelian group $T$ of a compact Lie group $K$, Vogan defined the corresponding finite root datum analogous to the classical root datum in \cite{hv}, and conjectured that the corresponding Weyl group of $T$ with respect to the identity component group of $K$ is generated by the set of root transvections in it. The conjecture remains open by now. For any compact simple Lie algebra $L_0$, Yu in \cite{y} classified all the abelian subgroups $T$ of the compact Lie group $\Aut(L_0)$ such that $T$ has the same dimension as that of its centralizer in $\Aut(L_0)$, which include all the maximal abelian subgroups of $\Aut(L_0)$, and also calculated all the Weyl groups.

For a complex simple Lie algebra $L$, the group $\Aut(L)$ is a simple algebraic group. One knows that there is a one-to-one correspondence between complex reductive algebraic groups and compact Lie groups. A $G$-grading $\Gamma$ on a simple Lie algebra $L$ is a fine grading (resp. a fine inner grading) if and only if the dual group $\wh{G}$ embeds as a maximal quasitorus in $\Aut(L)$ (resp. in $\Int(L)$). Thus the works in \cite{hv} and \cite{y} are closely related to the works mentioned in \textbf{1.1}. Furthermore, by (\ref{l1}) the Weyl group of the fine grading $(\Ga,G)$ and the Weyl group of the maximal quasitorus $\wh{G}$ are the same.\bigskip

\noindent \textbf{1.4.} A root system $R$ is a finite set of elements in a Euclidean space $V$ satisfying certain axioms. There is a one-to-one correspondence between reduced and irreducible root systems and simple complex Lie algebras (up to isomorphism).

In \cite{a}, for a Jordan subgroup $T$ of $\Aut(L)$ with $L$ a complex simple Lie algebra, an alternating bicharacter on $T$ is defined to study the Weyl group of $T$. In \cite{hy}, similar alternating bicharacters are also defined on elementary abelian 2-groups of compact simple Lie groups. An abelian group equipped with an alternating bicharacter is called a symplectic abelian group. Motivated by these works, in Section 3 we define a finite root system $(G,\bt,R)$ to be a finite set $R$ of elements in some finite symplectic abelian group $(G,\bt)$ satisfying certain axioms (Definition \ref{bb3}). Given any finite root system $(G,\bt,R)$, one always obtains a semisimple Lie algebra $L(R)$ with a standard (quasi-good) grading on it as well as an embedding $\wh{G}\hookrightarrow \Aut(L(R))$ (See Proposition \ref{l2}). The  finite root system $(G,\bt,R)$ is called \textit{good} if $\wh{G}$ is a finite maximal quasitorus in $\Aut(L(R))$ or in $\Int(L(R))$. Based on the main results in \cite{y}, in Theorem \ref{b1} we classify good finite maximal quasitorus $T$ in $\Aut(L)$ for any complex simple Lie algebra $L$. In Section 5 we show case by case that for any such pair $(L,T)$ one can always find some unique irreducible finite root system $(G,\bt,R)$, where $G=\wh{T}$ and $R$ is the set of roots of $T$ in $L$, such that $L\cong L(R)$ and the grading on $L$ induced by $T$ is isomorphic to the standard grading on $L(R)$. As a corollary, we classify all the reduced and irreducible finite root systems $R$ whose corresponding Lie algebra $L(R)$ is simple and the standard grading on it is a good inner grading.

It is known that Vogan's conjecture holds in these cases.  \bigskip

\noindent \textbf{1.5.} The paper is structured as follows. In Section 2 some results on twisted group algebras and symplectic abelian groups are collected, to prepare for Section 3, where we define finite root systems, corresponding quasi-good standings and Weyl groups. In Section 4 we will recall and prove some results about the gradings on simple Lie algebras. In Section 5 we give some typical examples of finite root systems $(G,\bt,R)$ and associated gradings, also give the embedding of $\wh{G}$ in $\Aut(L(R))$.  In Section 6 we classify good finite maximal quasitorus in $\Aut(L)$ for any complex simple Lie algebra $L$, and prove our main result Theorem \ref{b1}.  The proofs in the last section heavily rely on \cite{y}.

We use $\mathbb{Z}_n$ to denote $\mathbb{Z}/n\mathbb{Z}$ for any positive integer $n$. The other notations in the paper are standard.

\section{Twisted group algebras and symplectic abelian groups}
 \setcounter{equation}{0}\setcounter{theorem}{0}

 Let $G$ be a finite abelian group. Let $$A=\oplus_{\a\in G} \mathbb{C} u_a$$ be a vector space with basis $u_a$, and assume that there is a map $\xi:G\times G\rt \mathbb{C}^\times$ such that \be\lb{a1} u_a \cdot u_b=\xi(a,b) u_{a+b}, for\ all\ a,b\in G.\ee The following result can be found in Section 1, Chapter 2 of \cite{ka2}.

 \begin{lem}\lb{m6}
 $(A,\cdot)$ is a $G$-graded associative algebra if and only if   $$\xi(a,b)\xi(a+b,c)=\xi(a,b+c)\xi(b,c),for\ all\ a,b,c\in G,$$ i.e., $\xi\in Z^2(G,\mathbb{C}^\times).$ In this case $A$ is called a twisted group algebra of $G$ over $\mathbb{C}$.

 Assume that $\xi,\xi^{'}\in Z^2(G,\mathbb{C}^\times).$ Then
 $\xi$ and $\xi^{'}$ define isomorphic $G$-graded associative algebras if and only if $$\xi^{'}(a,b)/\xi(a,b)=\eta(a)\eta(b)/\eta(a+b)$$ for some $\eta:G\rt \mathbb{C}^\times$, i.e., $\xi$ and $\xi^{'}$ are cohomologous.

Thus the equivalence classes of twisted group algebras of $G$ over $\mathbb{C}$ are in 1-1 correspondence with $H^2(G,\mathbb{C}^\times)$.
  \end{lem}

A map $\bt:G\times G\rt  \mathbb{C}^\times$ is called an \textit{alternating bicharacter} on $G$ if  it is multiplicative in each variable and has the property that $\bt(a,a)=1$ for any $a\in G$, and is called a \textit{symmetric bicharacter} on $G$ if  it is multiplicative in each variable and has the property that $\bt(a,b)=\bt(b,a)$ for any $a,b\in G$. It is clear that all the  alternating bicharacters (resp. all the  symmetric bicharacters) on $G$ form an abelian group and will be denoted by $\w^2(G,\mathbb{C}^\times)$ (resp. $\textrm{Sym}^2(G,\mathbb{C}^\times)$). Note that $\w^2(G,\mathbb{C}^\times)$ is denoted by $P_{as}(G)$ in Chapter 8 of \cite{ka3}.

By Proposition 2.1 in Chapter 8 of \cite{ka3}, there is a short exact sequence of abelian groups

\bee \CD 1\rt \textrm{Ext}(G,\mathbb{C}^\times)\rt H^2(G,\mathbb{C}^\times)@>\psi>> \w^2(G,\mathbb{C}^\times)\rt 1,\endCD\eee
 where $\textrm{Ext}(G,\mathbb{C}^\times)=\{\overline{\al}\in H^2(G,\mathbb{C}^\times)|\al(a,b)=\al(b,a)\ for\ all\ a,b\in G\}$ which
 is 0 as $\mathbb{C}$ is algebraically closed and $\psi(\overline{\xi})(a,b)=\xi(a,b)\xi(b,a)^{-1}$. As $H^2(G,\mathbb{C}^\times)\cong Z^2(G,\mathbb{C}^\times)/B^2(G,\mathbb{C}^\times)$, the isomorphism $\psi$ can be written as a short exact sequence of abelian groups

 \be\lb{c5}\CD  1\rt B^2(G,\mathbb{C}^\times)\rt Z^2(G,\mathbb{C}^\times)@>\Psi>>  \w^2(G,\mathbb{C}^\times)\rt 1,\endCD\ee where the map $\Psi$ is \be \lb{h}\Psi: Z^2(G,\mathbb{C}^\times)\rt \wedge^2(G,\mathbb{C}^\times), \Psi(\xi)(a, b)= \xi(a,b) \xi(b,a)^{-1}.\ee

As $G$ is abelian, elements in $\Hom_\mathbb{Z}(G\ot G,\mathbb{C}^\times)$ are bimultiplicative maps $G\times G\rt \mathbb{C}^\times$, which are clearly 2-cocycles in $Z^2(G, \mathbb{C}^\times)$. Then (\ref{h}) restricts to a surjective group homomorphism
\be\lb{j} \Psi:\Hom_\mathbb{Z}(G\ot G,\mathbb{C}^\times)\rt  \w^2(G,\mathbb{C}^\times)\ee with the kernel
 consisting of all the symmetric bicharacters on $G$. Denote this kernel  by $ \textrm{Sym}^2(G,\mathbb{C}^\times)$. Thus one has \be  1\rt \textrm{Sym}^2(G,\mathbb{C}^\times)\rt\Hom_\mathbb{Z}(G\ot G,\mathbb{C}^\times)\rt \wedge^2(G,\mathbb{C}^\times)\rt 1.\ee

 If $\bt:G\times G\rt \mathbb{C}^\times$ is an alternating bicharacter on $G$, then $(G,\bt)$ is called a \textit{symplectic abelian group}. The radical of $\bt$ is
 $$\Rad(\bt)=\{a\in G|\bt(a,b)=1,\ for \ all\ b\in G\}.$$ One says that $\bt$ is \textit{nonsingular} if $\Rad(\bt)=0$.

By Lemma 2.7 and Corollary 2.10 in Chapter 8 of \cite{ka3}, one gets
\begin{prop}
One has $\mathbb{C}^\xi[G]\cong \oplus_{i=1}^k M(n,\mathbb{C})$, $k=|\Rad(\bt)|$ and $|G|=k n^2$. $\mathbb{C}^\xi[G]$ is simple if and only if $\bt$ is nonsingular.
\end{prop}
As a corollary one has
\begin{coro}\lb{n1}
For any $\beta\in\wedge^2(G,\mathbb{C}^\times)$, choose $\overline{\xi}\in H^2(G,\mathbb{C}^\times)$ with $\Psi({\xi})=\beta$, then we get a $G$-graded Lie algebra $L(\xi)$ obtained from the associative algebra $\mathbb{C}^\xi[G]$, i.e., $$[u_a,u_b]=(\xi(a,b)-\xi(b,a))u_{a+b}.$$ The Lie algebra structure of $L(\xi)$ depends only on $\beta$ and not on the $\xi$ chosen. One has $L\cong \gl(n,\mathbb{C})^{\op k}$ with $k=|\Rad(\bt)|$ and $|G|=k n^2$. In particular Isomorphic alternating bicharacters on $G$ corresponds to isomorphic $L(\xi)$.
  \end{coro}

 Let $\overline{\xi}\in H^2(G,\mathbb{C}^\times)$ and $\beta=\Psi(\xi)\in\wedge^2(G,\mathbb{C}^\times)$, then for any $a,b\in G$, $$\bt(a,b)=\xi(a,b)\xi(b,a)^{-1}.$$

 Let $L(\xi)$ be the reductive Lie algebra obtained from $\mathbb{C}^\xi[G]$. In $L(\xi)$ one has

 \begin{lem}\lb{i} (1) For any $a,b\in G$, $[u_a,u_b]=0$ if and only if $\xi(a,b)=\xi(b,a)$ if and only if $\bt(a,b)=1$. In particular, $[u_{m a},u_{ka}]=0$ for any $a\in G$ and $m,n\in \mathbb{Z}$.

(2) $L$ is commutative if and only if $\xi(a,b)=\xi(b,a) for\ all\ a,b\in G$ if and only if $\bt(a,b)\equiv 1 for\ all\ a,b\in G$.

(3) The center of $L$ is \be Z(L)=\oplus_{a\in \Rad(\bt)}\mathbb{C} u_a.\ee

(4) If $gcd(ord(a),ord(b))=1$ for some $a,b\in G$, then
 $\bt(a,b)=1$, $\xi(a,b)=\xi(b,a)$ and $[u_a,u_b]=0$.

  \end{lem}
\bp
(1) is obvious and (2),(3) follows from (1).

Let us prove (4). Assume $ord(a)=m$, $ord(b)=n$. Then $\bt(a,b)^m=\bt(ma,b)=1$ and $\bt(a,b)^n=\bt(a,nb)=1$. As $gcd(m,n)=1$, $\bt(a,b)=1$.  Then $\xi(a,b)=\xi(b,a)$ and $[u_a,u_b]=0$ by definition.
 \ep

Let $L^{'}=[L,L]$ be the derived Lie algebra of $L=L(\xi)$.

 \begin{prop}\lb{b}
               (1) $L^{'}=\oplus_{a\notin \Rad(\bt)}\mathbb{C} u_a$ is the semisimple ideal of $L$ and $L=Z(L)\oplus L^{'}$.

               (2) If $R$ is a subset of $G$ satisfying

               a) $R\st G\setminus \Rad(\bt)$ and $R$ generates $G$,

               b) If $a\in R$ then $-a\in R$, and

               c) If $\bt(a,b)\neq 1$ then $a+b\in R$,

               then $L(R)=\oplus_{a\in R}\mathbb{C} u_a$ is a semisimple Lie subalgebra of $L$.
\end{prop}

 \bp (1)Let $(,)$ be the Killing form on $L$. One knows by Corollary \ref{n1} that $L\cong\gl(n,\mathbb{C})^{\oplus k}$ for some integer $n,k$, and is reductive.

  If $a,b\notin \Rad(\bt)$ while $a+b\in \Rad(\bt)$, $\bt(a,b)=\bt(a,a+b)=1$ so $[u_a,u_b]=0$. Thus $L_1=\oplus_{a\notin \Rad(\bt)}\mathbb{C} u_a$ is a Lie subalgebra of $L$, and by Lemma \ref{i}.3, $$L=Z(L)\oplus L_1.$$ Then $L_1$ is an ideal of $L$ as $[L,L_1]=[L_1,L_1]\st L_1$. As $L^{'}=[L,L]=[L_1,L_1]\st L_1$, and $L=Z(L)\op L^{'}$, one must have $L^{'}=L_1$. So  $L^{'}=\oplus_{a\notin \Rad(\bt)}\mathbb{C} u_a$ is the semisimple ideal of $L$ and $L=Z(L)\oplus L^{'}$.

 (2) Let $(,)$ be the Killing form on $L(R)$. Assume $a,b\in R$. If $\bt(a,b)=1$ then $[u_a,u_b]=0$ by Lemma \ref{i}.1. If $\bt(a,b)\neq 1$ then $[u_a,u_b]=(\xi(a,b)-\xi(b,a))u_{a+b}\in L(R)$. So $L(R)$ is a Lie subalgebra of $L$.

It is clear that \bee (u_a,u_b)=0, if\ a+b\neq 0.\eee
 Choose $\xi\in \Hom_\mathbb{Z}(G\ot G,\mathbb{C}^\times)$ with $\Psi(\xi)=\bt$, i.e., $\bt(a,b)=\xi(a,b)\xi(b,a)^{-1}.$ Note that $\xi(a,b)$ is always a root of unity.
 One has  \be\lb{o}\begin{split}(u_a,u_{-a})&=tr(ad_{u_{-a}}\cdot ad_{u_a} )\\&=\sum_{b\in R} (\xi(a,b)-\xi(b,a))(\xi(-a,a+b)-\xi(a+b,-a))\\&=\sum_{b\in R} (2-\xi(a,b)\xi(b,-a)-\xi(b,a)\xi(-a,b)).\end{split}\ee  As $\xi(b,a)\xi(-a,b)=[\xi(a,b)\xi(b,-a)]^{-1}$ and $|\xi(b,a)\xi(-a,b)|=1$, one has
\[ 2-(\xi(a,b)\xi(b,-a)+\xi(b,a)\xi(-a,b))\geq 0,\] with equality holds if and only if  $\xi(a,b)\xi(b,-a)=1$, i.e., $\bt(a,b)=1$.

So $(u_a,u_{-a})=0$ if and only if $\bt(a,b)=1$ for any $b\in R$, which is equivalent to  $a\in \Rad(\bt)$ since $R$ generates $G$ by a). But $R\cap \Rad(\bt)=\varnothing$, so $(u_a,u_{-a})>0$ for any $a\in R$. Assume $x=\sum_{a\in R} \lambda_a u_a\ (\lambda_a\in \mathbb{C})$ is in the radical of the Killing form $(,)$, then for any $a\in R$, $(x,u_{-a})=\lambda_a(u_a,u_{-a})=0$ thus $\lambda_a=0$. So $x=0$ and the radical of the Killing form $(,)$ is 0. Therefore $(,)$ is nondegenerate and $L(R)$ is semisimple.
 \ep

\section{Finite root systems and corresponding quasi-good gradings on semisimple Lie algebras}
 \setcounter{equation}{0}\setcounter{theorem}{0}

If $T$ is a complex quasitorus, then one knows that $\wh{T}$ is a finitely generated abelian group, and the dual group of $\wh{T}$ is naturally isomorphic to $T$ by Pontryagin's duality. An algebraic group homomorphism of quasitori $\phi:T\rt S$ will induce a group homomorphism of finitely generated abelian groups $\wh{\phi}:\wh{S}\rt \wh{T}$, the dual homomorphism of which is just the algebraic group homomorphism $\phi:T\rt S$ under the canonical identification of a quasitorus with its bidual. The following result is clear.
\begin{lem}
An algebraic group homomorphism of quasitori $T\rt S$ is injective if and only if $\wh{S}\rt \wh{T}$ is surjective.
\end{lem}
\begin{lem}\lb{m8}
Assume that $L$ is a semisimple Lie algebra.

(1)An inclusion of quasitori $T\rt S$ of $\Aut(L)$ induces a surjective group homomorphism $\wh{S}\rt \wh{T}$ and a corresponding grading homomorphism.

(2) Conversely, a grading homomorphism on $L$ induced by the identity map $id:(\Ga_1,G)\rt (\Ga_2,H)$ induces an inclusion of quasitori $\wh{H}\rt\wh{G}$ of $\Aut(L)$.

 \end{lem}
 \bp
(1)The  inclusion of quasitori $T\rt S$ of $\Aut(L)$ clearly induces a group homomorphism $p:\wh{S}\rt \wh{T}$, which is surjective by last lemma. Let
 $$\Ga_1:L=\op_{g\in \wh{T}}L_g $$ and $$\Ga_2:L=\op_{g\in \wh{S}}L_g $$ be the respective gradings by $\wh{T}$ and by $\wh{S}$. If $g\in \wh{S}$ and $X\in L_g$, then it is easy to see that $X\in L_{p(g)}$ where $L_{p(g)}$ is the homogeneous component of degree $p(g)$ for $(\Gamma_2,\wh{T})$. Thus the identity map on $L$ induces a grading homomorphism.

(2) Conversely, a grading homomorphism induced by the identity map $id:(\Ga_1,G)\rt (\Ga_2,H)$ induces a surjective group homomorphism $G\rt H$ by definition, which induces an injective algebraic group homomorphism of quasitori $\wh{H}\rt\wh{G}$ by last lemma.
\ep

\begin{lem}\lb{x1} For a quasi-good $G$-grading $\Ga$ on $L$ induced by a finite quasitorus $T\st \Aut(L)$, there is a unique good $G_1$-grading $\Ga_1$ on $L$ induced by a finite maximal quasitorus $T_1\st \Aut(L)$ such that the identity map on $L$ induces a grading homomorphism $id:(\Ga_1,G_1)\rt (\Ga,G)$, which is injective on the support of $\Ga_1$. In particular, $T_1$ is the centralizer of $T$ in $\Aut(L)$, and $T\st T_1$.
\end{lem}
\bp Let $T_1$ be the centralizer of $T$ in $\Aut(L)$. Then
$T_1$ consists of semisimple automorphisms of $L$ preserving each 1-dimensional homogeneous component of $\Ga$, thus is a maximal quasitorus of $\Aut(L)$. Let $G_1=\wh{T_1}$. Then the $G_1$-grading $\Ga_1$ is good and has the desired property.

Next we will prove uniqueness. If  there is another good $G_2$-grading $\Ga_2$ on $L$ induced by a finite qusitorus $T_2\st \Aut(L)$ such that the identity map on $L$ induces a grading homomorphism $id:(\Ga_2,G_2)\rt (\Ga,G)$, then by Lemma \ref{m8} (2) $T\rt\ T_2$ is an inclusion of quasitorus. Then  $T_2\st Z(T)=T_1$ and then $T_2=T_1$ as  $T_2$ is a maximal quasitorus in $\Aut(L)$. The assertion $T\st T_1$ is clear.
\ep

Now we will define finite root systems.

 Let $G$ be a finite (additive) abelian group with an alternating bicharacter $$\bt:G\times G\rt \mathbb{C}^\times.$$ Then $(G,\bt)$ is a symplectic abelian group. Assume $\bt$ is nontrivial, i.e., $\Rad(\bt)\neq G$.
\begin{defi}\lb{bb3}
A subset $R$ of $G$ is a \textbf{finite root system} in $(G,\bt)$ if it satisfies
\begin{itemize}
  \item \textrm{FRS0}. $R\st G\setminus \Rad(\bt)$ and $R$ generates $G$.
 \item \textrm{FRS1}. If $a\in R$, then $-a\in R$.
  \item \textrm{FRS2}. If $\bt(a,b)\neq 1$, then $a+b\in R$.

\end{itemize}
  \end{defi}

We also say that $(G,\bt,R)$ is a finite root system. Let $Sp(G,\bt)$ denote the group of isometries of $(G,\bt)$.

The finite root system $R$ is called \textit{reduced} if the alternating bicharacter $\bt$ is nonsingular.

  Two finite root system $(G_i,\bt_i,R_i)$ ($i=1,2$) are \textit{isomorphic}, denoted $R_1\cong R_2$ or $(G_1,\bt_1,R_1)\cong (G_2,\bt_2,R_2)$, if there is a group isomorphism $\varphi:G_1\rt G_2$ preserving the respective alternating bicharacters and $\varphi(R_1)=R_2$.

 For a given finite abelian group $G$, if it admits a nonsingular alternating bicharacter then $G\cong H\times H$ for some abelian group $H$, and if there are two nonsingular alternating bicharacters $\bt_i $ on $G$ for $i=1,2$, then $(G,\bt_1)$ and $(G,\bt_2)$ are isometric. See Lemma 1.6 and Theorem 1.8 of \cite{ka3} for details. If $(G,\bt)$ is a nonsingular symplectic abelian group, then $R=G\setminus \{0\}$ is a finite root system in $(G,\bt)$ and the following result is clear.
\begin{lem}\lb{x5} If $(G_i,\bt_i,R_i)$ are reduced finite root systems with $R_i=G_i\setminus \{0\}$ for $i=1,2$. If $G_1\cong G_2$ then $(G_1,\bt_1,R_1)\cong (G_2,\bt_2,R_2)$.\end{lem}

Let $\varepsilon\in \mathbb{C}$ be a primitive $n$-th root of unity. Let $\mu_n=\{\varepsilon^i|i=0,\cdots,n-1\}$ be the cyclic subgroup of $\mathbb{C}^\times$ of order $n$. For any subgroup $H$ of $G$ of order $n$, define an action of $\mu_n$ on $H$ as follows:

\[ \varepsilon^i\cdot h=ih,\ \forall h\in H.\]

Given a finite root system $(G,\bt,R)$, assume $a\in R$ has order $n$. Fix a primitive $n$-th root of unity $\varepsilon$ and define $$s_a:G\rt G,b\mapsto b-\bt(a,b).a,$$ which is called a \textit{transvection} on $G$. It is directly verified that $s_a^n=1$ thus $s_a$ is invertible, and by (FRS2) one has $$s_a(R)=R,\forall a\in R.$$

\begin{lem} For any $a, b,c\in G$,
One has $$\bt(s_a(b),s_a(c))=\bt(b,c),$$ so $s_a\in Sp(G,\bt)$.
\end{lem}
\bp
Assume $a$ has order $n$ and $\varepsilon$ is a fixed primitive $n$-th root of unity. Then $\bt(a,b)=\ve^i,\ \bt(a,c)=\ve^j$ for some $i,j\in\mathbb{Z}$.\bee \begin{split} \bt(s_a(b),s_a(c))&=\bt(b-\bt(a,b).a, c-\bt(a,c).a)\\ &=\bt(b,c) \bt(a,c)^{-i}\bt(b,a)^{-j}\\ &=\bt(b,c) (\ve^j)^{-i}(\ve^{-i})^{-j}\\ &= \bt(b,c) \end{split}\eee
\ep

The \textit{Weyl group}  $W(R)$ of $R$ is defined to be the subgroup of $Sp(G,\bt)$  generated by $\{s_a|a\in R\}$.

 Let $L(R)=\oplus_{a\in R}\mathbb{C} u_a$ be the Lie algebra with the Lie bracket defined by some $\xi\in \Hom_\mathbb{Z}(G\ot G,\mathbb{C}^\times)$ such that $\Psi(\xi)=\bt$. By Proposition \ref{b} (2), $L(R)$ is a semisimple Lie algebra.

 The grading \[\Gamma:L(R)=\oplus_{a\in R}\mathbb{C} u_a\] is a $G$-grading on $L(R)$,$0\notin R$, and each of its homogeneous component has dimension 1. Thus $\Ga$ is a quasi-good grading on $L(R)$. We refer to $\Gamma$ as the \textit{standard} $G$-grading on $L(R)$ or the \textit{standard} $G$-grading of $R$. Together with Lemma \ref{x1} one has
 \begin{prop}\lb{l2}
 For any finite root system $(G,\bt,R)$, there corresponds to it a semisimple Lie algebra $L(R)$ with a standard quasi-good $G$-grading on it.

This grading is good if and only if  $\wh{G}$ is a maximal quasitorus in $\Aut(L(R))$, i.e., the centralizer of $\wh{G}$ in $\Aut(L(R))$ equals $\wh{G}$.
 \end{prop}

 A finite root system $(G,\bt,R)$  is called \textit{good} if $L=L(R)$ is simple and $\wh{G}$ is a  finite maximal quasitorus in $\Aut(L)$ or in $\Int(L)$.
 A finite root system $R$ in $(G,\bt)$ is called \textit{reducible} if $G$ is an orthogonal product of two subgroups $G_1$ and $G_2$, $R$ is a disjoint union of two  nonempty orthogonal subsets $R_1$ and $R_2$, and $R_i$ is a finite root system in $G_i$ for $i=1,2$. It is clear that in this case $L(R)=L(R_1)\op L(R_2)$ is a direct sum of semisimple ideals $L(R_1)$ and $L(R_2)$.  A finite root system $R$ in $(G,\bt)$ is \textit{irreducible} if it is not reducible. In particular, if $(G,\bt,R)$  is good then $L(R)$ is simple thus it must be irreducible.

 \begin{lem}\lb{l3}
Assume that $(G,\bt,R)$ is a finite root system and $H$ is any subgroup of $\Rad(\bt)$. Let $\overline{G}=G/H$, $\overline{R}$ be the image of $R$ in $\overline{G}$, and $\overline{\bt}$ be the alternating bicharacter on $\overline{G}$ induced by $\bt$, i.e., $\overline{\bt}(\overline{g},\overline{h})=\bt(g,h)$ for any $\overline{g},\overline{h}\in \overline{G}$.

(1) One has that $(\overline{G},\overline{\bt},\overline{R})$ is a finite root system. If $H=\Rad(\bt)$ then $\overline{R}$ is a reduced finite root system, and is referred as the corresponding reduced finite root system of $R$.

(2) There is a surjective Lie algebra homomorphism from $L(R)$ to $L(\overline{R})$ that maps $L(R)_g$ to $L(\overline{R})_{\overline{g}}$ for any $g\in R$. If $L(R)$ is simple then this is an isomorphism of Lie algebras, and in this case the homogeneous components of $L(R)$ and $L(\overline{R})$ are the same.
  \end{lem}
\bp (1)  is clear.

(2) Let $\overline{\xi}\in Hom_\mathbb{Z}(\overline{G}\ot \overline{G}, \mathbb{C}^\times)$ such that $\Psi(\overline{\xi})=\overline{\bt}$. Let  ${\xi}\in Hom_\mathbb{Z}({G}\ot {G}, \mathbb{C}^\times)$ be the pull-back of $\overline{\xi}$, i.e., for any $g,h\in G$, $\xi(g,h)=\overline{\xi}(\overline{g},\overline{h})$. Then it is directly verified that $\Psi({\xi})={\bt}$.

Then $\phi:L(R)\rt L(\overline{R}), u_g\mapsto u_{\overline{g}}$ defines a surjective Lie algebra homomorphism, as $$[u_g,u_h]=(\xi(g,h)-\xi(h,g))u_{g+h}$$ and $$[u_{\overline{g}},u_{\overline{h}}]=(\overline{\xi}(\overline{g},\overline{h})-\overline{\xi}(\overline{h},\overline{g}))u_{\overline{g}+\overline{h}}
=(\xi(g,h)-\xi(h,g))u_{\overline{g+h}}.$$

The last assertion in (2) is clear.
\ep

If two finite root systems $(G_i,\bt_i,R_i)$ are isomorphic for $i=1,2$, then it is clear that $G_1\cong G_2$ and the corresponding standard gradings are isomorphic. The converse may not hold in general, but it holds in a special case.

\begin{lem}\lb{x7}
Assume that the standard gradings of two finite root systems $(G_i,\bt_i,R_i)$   are isomorphic for $i=1,2$ and $G_1\cong G_2\cong \mathbb{Z}_2^n$ for some positive integer $n$. Then $(G_1,\bt_1,R_1)\cong (G_2,\bt_2,R_2)$.
\end{lem}
\bp
By assumption, there is a Lie algebra isomorphism $\psi:L(R_1)\rt L(R_2)$ and a group isomorphism $\phi:G_1\rt G_2$ such that $\phi(R_1)=R_2$ and for any $a\in R_1$, $\psi(L(R_1)_a)=L(R_2)_{\phi(a)}$, i.e., $\psi(u_a)=\lambda_a u_{\phi(a)}$ for some $0\neq \lambda_a\in \mathbb{C}$. As $\psi$ preserves Lie brackets, for any $a,b\in R_1$,
$[u_a,u_b]=0$ if and only if $[u_{\phi(a)},u_{\phi(b)}]=0$. Then by Lemma \ref{i} (1), for any $a,b\in R_1$, \be\lb{x3}\bt_1(a,b)=1\ if\ and\ only\ if\ \bt_2({\phi(a)},{\phi(b)})=1.\ee As $G_1\cong G_2\cong \mathbb{Z}_2^n$, any nonzero element in $G_1$ and $G_2$ has order 2, thus $\bt_1$ and $\bt_2$ take values in $\pm 1\in \mathbb{C}^\times.$
So (\ref{x3}) implies  for any $a,b\in R_1$, \[\lb{x2}\bt_1(a,b)=-1\ if\ and\ only\ if\ \bt_2({\phi(a)},{\phi(b)})=-1.\] As $R_1$ generates $G_1$, one has $\bt_2({\phi(a)},{\phi(b)})=\bt_1(a,b)$ for any $a,b\in G_1$. Therefore  $\phi:G_1\rt G_2$ preserves the respective alternating bicharacters and $\phi(R_1)=R_2$ thus $\phi$ is an isomorphism of the two finite root systems $(G_i,\bt_i,R_i)$ for $i=1,2$.
\ep

Finally we will give two more examples of finite root systems.

 Let $q:\mathbb{F}_2^{2k}\rt \mathbb{F}_2$ be a  quadratic form and \[\al_1:\mathbb{F}_2^{2k}\times \mathbb{F}_2^{2k}\rt \mathbb{F}_2, \al_1(a,b)=q(a+b)-q(a)-q(b)\] be its \textit{polarization}. Assume that $q$ is nonsingular, i.e.  $\al_1$ is a nonsingular alternating bilinear form.  One knows that there are two types of nonsingular quadratic forms on $\mathbb{F}_2^{2k}$ up to isomorphism. We will always identify the group $G=\mathbb{Z}_2^{2k}$ with the additive group of $\mathbb{F}_2^{2k}$.

\begin{lem}\lb{m4} Assume that $\al_1$ is nonsingular. Let \[\al: G\times G\rt \mathbb{C}^\times, \al(a,b)=(-1)^{\al_1(a,b)}. \] Assume that \[R=\{a\in G\setminus \{0\}|q(a)=1\}\] generate $G$. Then $R$ is a finite root system in the symplectic abelian group $(G,\al)$.
\end{lem}
\bp
 FRS0 and FRS1 are obvious. Now we verify FRS2. Assume $a,b\in R$ and $\al(a,b)\neq 1$. Then $\al(a,b)=-1$, $\al_1(a,b)=1$ and $q(a+b)=q(a)+q(b)+\al_1(a,b)=1$, so $a+b\in R$.
\ep

Let \[
\bt_1:\mathbb{F}_2^{2k}\times \mathbb{F}_2^{2k}\rightarrow \mathbb{F}_2,\bt_1(a,b)={\sum_{i=1}^k(a_{2i}b_{2i-1}-a_{2i-1}b_{2i})}
\] be the unique nonsingular alternating bilinear form on $\mathbb{F}_2^{2k}$ up to isomorphism. Let \[
\bt:G\times G\rightarrow \mathbb{C}^{\times},\bt(a,b)=(-1)^{\bt_1(a,b)}.
\] Then $(G,\bt)$ is a nonsingular symplectic abelian group. Up to isomorphism there are two quadratic forms on $\mathbb{F}_2^{2k}$ that polarize to $\bt_1$,  \be\lb{aa8} g:\mathbb{F}_2^{2k}\rightarrow \mathbb{F}_2,g(a)={\sum_{i=1}^{k}a_{2i-1}a_{2i}},\ee and  \be\lb{aa9} f:\mathbb{F}_2^{2k}\rightarrow \mathbb{F}_2,f(a)={a_1^2+a_2^2+\sum_{i=1}^{k}a_{2i-1}a_{2i}}.\ee So we have the following two examples of finite root systems and corresponding standard gradings.

\begin{lem}\lb{b9}
Assume $G=\mathbb{Z}_2^{2k}$ with $k\geq 3$. Let $R=\{a\in G|g(a)=1\}$.
Then  $R$ is a finite root system in $(G,\bt)$.
 \end{lem}

\bp The subset \[B=\{e_{2i-1}+e_{2i},\  e_{2i-1}+e_{2i}+e_{2i+1}|i=1,\cdots,k\}  \ (e_{2k+1}=e_1)\] of $R$ generates $G$, thus by last lemma  $R$ is a finite root system in $(G,\bt)$.

\ep

 One will see that
\[L(R)\cong \so(2^k,\mathbb{C})\] by Lemma \ref{a4} (1), and the corresponding grading is (\ref{z1}).

\begin{lem}\lb{bb8}
Assume  $G=\mathbb{Z}_2^{2k}$ with $k\geq 1$. Let $R=\{a\in G|f(a)=1\}$.Then  $R$ is a finite root system in $(G,\bt)$.
 \end{lem}
\bp It is directly verified that \bee B=\begin{cases} \{e_1\}\cup\{e_1+e_2\} & \text{if\ k=1};\\
\{e_1,e_{2i-1}+e_{2i}|i=1,\cdots,k\}\cup \{e_{2i-1}+e_{2i}+e_{2i+1}|i=1,\cdots,k-1\} & \text{if\ $k>1$.}
\end{cases} \eee is contained in $R$ and generates $G$, thus  $R$ is a finite root system in $(G,\bt)$.
\ep
 One will see that
\[L(R)\cong\sp(2^k,\mathbb{C})\] by Lemma \ref{a4} (2), and the corresponding grading is (\ref{z2}).

\section{Some results about the gradings on simple Lie algebras}
 \setcounter{equation}{0}\setcounter{theorem}{0}

The gradings on classical simple Lie algebras are all described in \cite{bsz}, which can be obtained by the gradings on matrix algebras.

 For a nonsingular symmetric or skew-symmetric bilinear form $\phi$ on $\mathbb{C}^n$, one knows that the adjoint map $*:M(n,\mathbb{C})\rt M(n,\mathbb{C})$ defined by
  $$\phi(Xu,v)=\phi(u,X^* v) $$ is an involution (i.e., involutive anti-automorphism) on $M(n,\mathbb{C})$. If  $\Phi$ is the matrix of $\phi$ with respect to the standard basis of $\mathbb{C}^n$, then in matrix form one has $$X^*=\Phi^{-1} X^t \Phi,$$
Conversely one can show that any involution on $M(n,\mathbb{C})$ can be defined in this way by some nonsingular symmetric or skew-symmetric bilinear form on $\mathbb{C}^n$. See Section 5 of \cite{bsz} for details. We call $(M(n,\mathbb{C}),*)$ an \textit{involutive} matrix algebra. In the remaining of the section, $M$ will be assumed to be a matrix algebra.

Assume that $(M,*)$ is an {involutive} matrix algebra. Let \[\Aut(M,*)=\{\si\in \Aut(M)|\si\circ *=*\circ\si\}.\] If $H$ is a subgroup of $\Aut(M,*)$, then its action on $M$ is compatible with $*$,i.e.,$h\circ *=*\circ h$ for any $h\in H$, and we say that $(M,*)$ has an $H$-action. If $M$ has a $G$-grading compatible with $*$,i.e., $M_g^*=M_g$ for any $g\in G$, then we say that $(M,*)$ has a $G$-grading. It is clear that $(M,*)$ has a $G$-grading if and only if $(M,*)$ has a $\wh{G}$-action.

Let \[K(M,*)=\{X\in M|X^*=-X\}\] and \[H(M,*)=\{X\in M|X^*=X\}.\] The subspace $K(M,*)$ is closed under the Lie bracket $[a,b]=ab-ba$ and is always regarded as a Lie algebra with this Lie bracket in the paper. A $G$-grading on $(M,*)$ clearly induces a $G$-grading on $K(M,*)$ (resp. on $H(M,*)$).

One knows that simple Lie algebras of type \textbf{B}, \textbf{C} and \textbf{D} can be realized as some $K(M,*)$ for suitable $(M,*)$, where $M=M(n,\mathbb{C})$ and $*$ corresponds to some nonsingular symmetric or skew-symmetric bilinear form $\phi$ on $\mathbb{C}^n$.

 Let $G=\mathbb{Z}_2\times \mathbb{Z}_2$,$M=M(2,\mathbb{C})$. Let
\be X_2=\left(
             \begin{array}{cc}
               -1 & 0\\
               0 & 1 \\
             \end{array}
           \right) \an \
           Y_2=\left(
             \begin{array}{cc}
               0 & 1\\
               1 & 0 \\
             \end{array}
             \right).
\ee
 Let $Z_{(i,j)}=X_2^i Y_2^j$ for $(i,j)\in G$.
The grading $M=\op_{a\in G} M_a$, where $M_a=\mathbb{C} Z_a$ for any $a\in G$, is referred as the standard $G$-grading on $M$. See  Example \ref{c2}.

Assume that $*$ is an involution on $M$ corresponding to the skew-symmetric form with the  matrix $\Phi_1=\left(
             \begin{array}{cc}
               0 & 1\\
               -1 & 0 \\
             \end{array}
           \right)$. By Lemma 3 (1) in \cite{bsz}, $(M,*)$ has the $G$-grading with \be \lb{a1}K(M,*)=\{\left(
             \begin{array}{cc}
               a & b\\
               c & -a \\
             \end{array}
           \right)\}=M_{(1,0)}\op M_{(0,1)}\op M_{(1,1)} \ee and \[H(M,*)=\{\left(
             \begin{array}{cc}
               a & 0\\
               0 & a \\
             \end{array}
           \right)\}=M_{(0,0)}.\]

 Assume that $*$ is an involution on $M$ corresponding to the symmetric form with the matrix $I_2=\left(
             \begin{array}{cc}
               1 & 0\\
               0 & 1 \\
             \end{array}
           \right)$.

By Lemma 3 (3) of \cite{bsz}, $(M,*)$ has the $G$-grading with \be \lb{a2}K(M,*)=\{\left(
             \begin{array}{cc}
               0 & b\\
               -b & 0 \\
             \end{array}
           \right)\}=M_{(1,1)}\ee and \[H(M,*)=\{\left(
             \begin{array}{cc}
               a & b\\
               b & c \\
             \end{array}
           \right)\}=M_{(0,0)}\op M_{(0,1)}\op M_{(1,0)}.\]
\bigskip

Next let us introduce the direct product gradings on associative algebras. Assume for $i=1,2$, $A_i$ is an associative algebra, and $(\Ga_i,G_i)$ is the respective gradings on $A_i$, then $(\Ga, G_1\times G_2)$ is the \textit{direct product grading} on $A=A_1\otimes A_2$ with $$\Gamma:A_1\ot A_2=\bigoplus_{(a,b)\in G_1\times G_2 } (A_1\otimes A_2)_{(a,b)}$$ where $$(A_1\otimes A_2)_{(a,b)}=(A_1)_a\otimes (A_2)_b, (a,b)\in G_1\times G_2.$$

If $(M_i,*_i)$ is an involutive matrix algebra for $i=1,2$, then their direct product $(M_1,*_1)\ot (M_2,*_2)$ is defined to be the involutive matrix algebra $(M_1\ot M_2,*_1\ot *_2)$, where $(C\ot B)^{*_1\ot *_2}=C^{*_1}\ot B^{*_2}$ for any $C\ot B\in M_1\ot M_2$. If $(M_i,*_i)$ has a $G_i$-grading, then the $G_1\times G_2$-grading on $M_1\ot M_2$ is compatible with $*_1\ot *_2$, thus $(M_1\ot M_2,*_1\ot *_2)$ has the naturally defined $G_1\times G_2$-grading.

Assume that $(M_i,*_i)$ has a $\mathbb{Z}_2\times \mathbb{Z}_2$-grading for $i=1,\cdots,k$, where $M_i=M(2,\mathbb{C})$ and the $\mathbb{Z}_2\times \mathbb{Z}_2$-grading on $M_i$ is standard for all $i$. Then the tensor product $(M,*)$ of $(M_i,*_i)$ for $i=1,\cdots,k$ has a $G$-grading with \[G=(\mathbb{Z}_2\times \mathbb{Z}_2)^k=\mathbb{Z}_2^{2k}.\] Thus $K(M,*)$ also has a $G$-grading. Assume that $*_i$ corresponds to a nonsingular bilinear form $\phi_i$, then $*$ corresponds to the nonsingular bilinear form $\phi=\phi_1\ot\cdots\ot\phi_k$.  Assume the total number of skew-symmetric factors in $\phi$ is $m$. If $m$ is even (resp. odd), then $\phi$ is symmetric (resp. skew-symmetric).  Next we will compute the support $R$ of the $G$-grading on $K(M,*)$ in the two cases: $m=0$ or $m=1$. (It is not hard to show that the grading on $K(M,*)$ is isomorphic to the first case if $m$ is even, and is isomorphic to the second case if $m$ is odd.)

Let \be\lb{o3}\bt_1: \mathbb{F}_2^{2k}\times  \mathbb{F}_2^{2k}\rightarrow \mathbb{F}_2,\bt_1(a,b)={\sum_{i=1}^k(a_{2i}b_{2i-1}-a_{2i-1}b_{2i})} \ee be a nonsingular bilinear form on $\mathbb{F}_2^{2k}$.

Recall the nonsingular quadratic forms $g$ and $f$ on $\mathbb{F}_2^{2k}$ defined in (\ref{aa8}) and (\ref{aa9}). One knows that they both polarize to $\bt_1$.

 Let  \be\lb{l5} \bt: \mathbb{F}_2^{2k}\times  \mathbb{F}_2^{2k}\rightarrow \mathbb{C}^{\times},\bt(a,b)=(-1)^{\bt_1(a,b)}=(-1)^{{\sum_{i=1}^k(a_{2i}b_{2i-1}-a_{2i-1}b_{2i})}}.\ee
\begin{lem}\lb{a4}
Let $k\geq 1$ be an integer and $G=\mathbb{Z}_2^{2k}$.

(1) Assume  the total number of skew-symmetric factors $\phi_i$ in $\phi$ is 0 and the matrix for each $\phi_i$ is $I_2$. Then $K(M,*)\cong \so(2^k,\mathbb{C})$, $R=\{a\in G\setminus \{0\}|g(a)=1\}$. One has $|R|=2^{2k-1}-2^{k-1}$, and this grading on $\so(2^k,\mathbb{C})$ is
\be \lb{z1} \so(2^k,\mathbb{C})=\bigoplus_{a\in G,g(a)=1}\mathbb{C} X_2^{a_1}Y_2^{a_2}\ot X_2^{a_3}Y_2^{a_4}\ot\cdots\ot X_2^{a_{2k-1}}Y_2^{a_{2k}}.\ee

Conversely, let $R$ be as above. Then $R$ is a finite root system in $(G,\bt)$ and the standard grading on $L(R)$ is isomorphic to the $G$-grading on $\so(2^k,\mathbb{C})$. Moreover, this grading is good if $k\geq 3$.

(2) Assume the total number of skew-symmetric factors $\phi_i$ in $\phi$ is 1, the matrix of $\phi_1$ is $\Phi_1$ and the matrix of $\phi_i$ is  $I_2$ for $i>1$. Then $K(M,*)\cong \sp(2^k,\mathbb{C})$, $R=\{a\in G\setminus \{0\}|f(a)=1\}$. One has $|R|=2^{2k-1}+2^{k-1}$, and this grading on $\sp(2^k,\mathbb{C})$ is
\be \lb{z2} \sp(2^k,\mathbb{C})=\bigoplus_{a\in G,f(a)=1}\mathbb{C} X_2^{a_1}Y_2^{a_2}\ot X_2^{a_3}Y_2^{a_4}\ot\cdots\ot X_2^{a_{2k-1}}Y_2^{a_{2k}}.\ee

Conversely, let $R$ be as above. Then $R$ is a finite root system in $(G,\bt)$ and the standard grading on $L(R)$ is isomorphic to the $G$-grading on $\sp(2^k,\mathbb{C})$. Moreover, this grading is good.
\end{lem}
\bp
First let us consider the case $k=1$. Then $G=\mathbb{Z}_2\times \mathbb{Z}_2$ and $M=M(2,\mathbb{C})$. Then by (\ref{a2}) the support of $\so(2,\mathbb{C})=K(M,*)$ is \be\lb{m2}\{(1,1)\}=\{(x_1,x_2)\in G|x_1x_2=1\}.\ee By (\ref{a1}) the support of $ \sp(2,\mathbb{C})=K(M,*)$ is \be\lb{m3}\{(0,1),(1,0),(1,1)\}=\{(x_1,x_2)\in G|x_1^2+x_2^2+x_1x_2=1\}.\ee

Next let us consider the general case. Then $G=\mathbb{Z}_2^{2k}$. Let $C=C_1\ot\cdots\ot C_k\in M(2,\mathbb{C})\ot\cdots\ot M(2,\mathbb{C})$, $C_i$ being one of $Z_a$ with $a\in \mathbb{Z}_2\times \mathbb{Z}_2$. Let $v_1\ot\cdots\ot v_k, u_1\ot\cdots\ot u_k\in (\mathbb{C}^2)^{\ot k}$. Then \be\lb{m1}\phi_i(C_i(v_i),u_i)=\begin{cases} -\phi_i(v_i,C_i(u_i))&\text{if $C_i\in K(M, *_i)$;}\\ \phi_i(v_i,C_i(u_i))&\text{if $C_i\in H(M, *_i)$.}\end{cases}\ee
One has
 \[\phi( C_1\ot\cdots\ot C_k(v_1\ot\cdots\ot v_k),u_1\ot\cdots\ot u_k)=\phi_1(C_1(v_1),u_1)\cdots \phi_k(C_k(v_k),u_k),\] and
\[\phi(v_1\ot\cdots\ot v_k, C_1\ot\cdots\ot C_k(u_1\ot\cdots\ot u_k))=\phi_1(v_1,C_1(u_1))\cdots \phi_k(v_k,C_k(u_k)).\]
Thus by (\ref{m1}) \[\phi( C_1\ot\cdots\ot C_k(v_1\ot\cdots\ot v_k),u_1\ot\cdots\ot u_k)=-\phi(v_1\ot\cdots\ot v_k, C_1\ot\cdots\ot C_k(u_1\ot\cdots\ot u_k))\]
if and only if there are exactly odd number of $C_i\in K(M, *_i)$.

(1) Assume  the total number of skew-symmetric factors $\phi_i$ in $\phi$ is 0 and the matrix for each $\phi_i$ is $I_2$.  Recall \[g:\mathbb{Z}_2^{2k}\rt \mathbb{Z}_2,g(x)=\sum_{i=1}^k x_{2i-1}x_{2i}.\] Let $a=(a_1,a_2,\cdots,a_{2k-1},a_{2k})\in \mathbb{Z}_2^{2k}$. Then by (\ref{m2}) $a\in R$ if and only if $g(a)=1$, so the support $R$ of $K(M,*)$ is those $a\in G$ with $g(a)=1$. Similarly the support of $H(M,*)$ is those $a\in G$ with $g(a)=0$. The fact  $|R|=2^{2k-1}-2^{k-1}$ is simple and we omit its proof. The grading of $K(M,*)\cong \so(2^k,\mathbb{C})$ is clearly (\ref{z1}).

Recall the definition of $\bt$ in (\ref{l5}). Then $\bt\in\w^2(G,\mathbb{C}^\times)$ is a nonsingular alternating bicharacter on $G$. By Lemma \ref{b9}, $R$ is a finite root system in $(G,\bt)$.
Let $\xi:\mathbb{Z}_2^{2k}\times \mathbb{Z}_2^{2k}\rt \mathbb{C}^\times, \xi(a,b)=(-1)^{\sum_{i=1}^k a_{2i}b_{2i-1}}.$ Then $\Psi(\xi)=\bt$.

Conversely, define $\phi:L(R)\rt \so(2^k,\mathbb{C})$, \[\phi(u_a)=Z_{(a_1,a_2)}\ot Z_{(a_3,a_4)}\ot\cdots\ot Z_{(a_{2k-1},a_{2k})}.\]
where $a=(a_1,a_2,\cdots,a_{2k-1},a_{2k}).$ Then it is easy to verify that this is an isomorphism of $G$-gradings as follows.

\be\lb{d1}\begin{split}[\phi(u_a),\phi(u_b)]&=[Z_{(a_1,a_2)}\ot Z_{(a_3,a_4)}\ot\cdots\ot Z_{(a_{2k-1},a_{2k})},Z_{(b_1,b_2)}\ot Z_{(b_3,b_4)}\ot\cdots\ot Z_{(b_{2k-1},b_{2k})}]
 \\&=Z_{(a_1,a_2)}Z_{(b_1,b_2)}\ot\cdots\ot Z_{(a_{2k-1},a_{2k})}Z_{(b_{2k-1},b_{2k})}-Z_{(b_1,b_2)}Z_{(a_1,a_2)}\ot\cdots\ot Z_{(b_{2k-1},b_{2k})}Z_{(a_{2k-1},a_{2k})}\\&=[\xi(a,b)-\xi(b,a)]Z_{(a_1+b_1,a_2+b_2)}\ot\cdots\ot Z_{(a_{2k-1}+b_{2k-1},a_{2k}+b_{2k})}\\&=[\xi(a,b)-\xi(b,a)]\phi(u_{a+b})\\&=\phi([u_a,u_b]).\end{split}\ee

Assume $k\geq 3$. As $PO(2^k,\mathbb{C})\st PGL(2^k,\mathbb{C})$, one knows that $\wh{G}\cong \mathbb{Z}_2^{2k}$ embeds as a maximal quasitorus in $\Aut(\so(2^k,\mathbb{C}))\cong PO(2k,\mathbb{C})$, so this grading is a good grading. (If $k=1$ or 2 then $L(R)$ is not simple.)

(2) Assume the total number of skew-symmetric factors $\phi_i$ in $\phi$ is 1, the matrix of $\phi_1$ is $\Phi_1$ and the matrix of $\phi_i$ is  $I_2$ for $i>1$. Recall \[f:\mathbb{Z}_2^{2k}\rt \mathbb{Z}_2,f(x)=(x_1^2+x_2^2+x_1 x_2)+\sum_{i=2}^k x_{2i-1}x_{2i}=x_1^2+x_2^2+\sum_{i=1}^k x_{2i-1}x_{2i}.\]  Then  by (\ref{m2}) and (\ref{m3}) $a=(a_1,a_2,\cdots,a_{2k-1},a_{2k})\in R$ if and only if $f(a)=1$, so the support $R$ of $K(M,*)$ is those $a\in G$ with $f(a)=1$. Similarly the support of $H(M,*)$ is those $a\in G$ with $f(a)=0$. The fact  $|R|=2^{2k-1}+2^{k-1}$ is simple and we omit its proof.  The grading of $K(M,*)\cong \so(2^k,\mathbb{C})$ is clearly (\ref{z2}).

Let $\bt$ be as above. By Lemma \ref{bb8}, $R$ is a finite root system in $(G,\bt)$. Let \[\xi^{'}:G\times G\rightarrow \mathbb{C}^{\times},\xi^{'}(a,b)=(-1)^{a_1b_1+a_2b_2+\sum_{i=1}^k a_{2i}b_{2i-1}},\] where $a=(a_1,\cdots,a_{2k}), b=(b_1,\cdots,b_{2k})\in G$. Then $\Psi(\xi^{'})=\bt$. By Lemma \ref{m6} and (\ref{c5}), there is some $\eta: G\rt \mathbb{C}^\times$ satisfying $\xi^{'}(a,b)/\xi(a,b)=\eta(a)\eta(b)/\eta(a+b)$ for any $a,b\in G$.

Conversely, define $\phi:L(R)\rt  \sp(2^k,\mathbb{C})$, \[\phi(u_a)=\eta{(a)}Z_{(a_1,a_2)}\ot Z_{(a_3,a_4)}\ot\cdots\ot Z_{(a_{2k-1},a_{2k})}.\]
where $a=(a_1,a_2,\cdots,a_{2k-1},a_{2k})\in R.$ Then one can verify that this is an isomorphism of $G$-gradings as in (\ref{d1}).

As $PSp(2^k,\mathbb{C})\st PGL(2^k,\mathbb{C})$, one knows that $\wh{G}\cong \mathbb{Z}_2^{2k}$ embeds as a maximal quasitorus in $\Aut( \sp(2^k,\mathbb{C}))\cong PSp(2^k,\mathbb{C})$, so this grading is good.

\ep

\section{Good finite maximal quasitorus in $\Aut(L)$ for simple Lie algebras $L$ and corresponding finite root systems  }
 \setcounter{equation}{0}\setcounter{theorem}{0}
Given a complex simple Lie algebra $L$, a finite maximal quasitorus $T$ of $\Aut(L)$ or $\Int(L)$ is said to be good if the induced grading on $L$ is quasi-good, i.e.,
$\dim L_{\alpha}\leq 1$ for every character $\alpha\in \wh{T}$. In this case, the set of roots $R$ of $T$ in $L$ is a subset in $G=\wh{T}$, and we will find some alternating bicharacter $\bt$ on $G$ such that $(G,\bt,R)$ is a finite root system, $L\cong L(R)$ and the grading on $L$ induced by $T$ is isomorphic to the standard grading on $L(R)$. It will be seen that the finite root system is unique up to isomorphism. It is a pity that we do not find a canonical way to construct the alternating bicharacter.

Recall that a finite root system $(G,\bt,R)$ is good if $L=L(R)$ is simple and $\wh{G}$ is a  finite maximal quasitorus in $\Aut(L)$ or in $\Int(L)$. In the following examples, we will construct good finite root systems $(G,\bt,R)$, and identify $L=L(R)$ with the corresponding standard gradings on it. Then we will give the embedding of $T=\wh{G}$ in $\Aut(L)$ or $\Int(L)$. We also give $W(R)$ in each case. People have known that  $W_K(\Ga)=W(R)$ always holds in these cases, where $K=\Int(L)$.

In Section 6 we will prove that these exhaust all the good finite root systems.\bigskip

\noindent \textbf{5.1.} Reduced finite root systems $(G,\bt,R)$ with $R=G\setminus\{0\}$

\begin{ex}\lb{c2}
 Let $G=\mathbb{Z}_n\times \mathbb{Z}_n$. Let $\varepsilon$ be a primitive $n$-th root of unity. For any $(i,j),(s,t)\in \mathbb{Z}_n\times \mathbb{Z}_n=G$, let $\bt((i,j),(s,t))=\varepsilon^{js-it}.$ It is clear that $\bt$ is a nonsingular alternating bicharacter on $G$ and $ R=G\setminus \{0\}$ is a finite root system in $(G,\bt)$.

 Let $\xi\in Hom_\mathbb{Z}(G\otimes G, \mathbb{C}^\times)$ be defined by $\xi((i,j),(s,t))=\varepsilon^{js}.$ Then $\Psi(\xi)=\bt$.

Next we will show $L(R)\cong \sl(n,\mathbb{C})$.

Let $M=M(n,\mathbb{C})$ and $G=\mathbb{Z}_n\times \mathbb{Z}_n$. Let \be \lb{y}X=X_n=\begin{pmatrix}
\varepsilon^{n-1} & 0 & \cdots & 0 \\
0 & \varepsilon^{n-2} & \cdots & 0 \\
 &  & \cdots &  \\
0 & 0 & \cdots & 1
\end{pmatrix},\quad Y=Y_n=\begin{pmatrix}
0 & 1 & 0 & \cdots & 0 \\
0 & 0 & 1 & \cdots & 0 \\
 &  &  & \cdots &  \\
0 & 0 & 0 & \cdots & 1 \\
1 & 0 & 0 & \cdots & 0
\end{pmatrix}. \ee
Then $XY=\varepsilon YX$ and $X^n=Y^n=I_n$ where $I_n$ will always denote the $n\times n$ identity matrix. Then \be\Ga: M=\op_{(i,j)\in G} \mathbb{C} X^i Y^{-j}\ee is a $G$-grading on $M$, where $M_{(i,j)}=\mathbb{C} X^i Y^{-j}$, called the $\varepsilon$-grading on $M$ in \cite{bsz}. 

 Then

\be \sl(n,\mathbb{C})=\op_{(i,j)\in G\setminus \{0\}} \mathbb{C} X^i Y^{-j}\ee is a $G$-grading on $\sl(n,\mathbb{C})$.

  Define \[\varphi:L(R)\rt \sl(n,\mathbb{C}), u_{i,j}\mapsto X^i Y^{-j}.\] \begin{align*}[\varphi(u_{i,j}),\varphi(u_{s,t})]=&X^iY^{-j}X^sY^{-t}-X^sY^{-t}X^iY^{-j}\\=&(\varepsilon^{js
}-\varepsilon^{ti})X^{i+s}Y^{-(j+t)}\\=&(\xi((i,j),(s,t))-\xi((s,t),(i,j)))X^{i+s}Y^{-(j+t)}\\=&(\xi((i,j),(s,t))-\xi((s,t),(i,j)))\varphi(u_{i+s,j+t})
\\=&\varphi([u_{i,j},u_{s,t}]).
\end{align*}
Thus $\varphi:L(R)\rt \sl(n,\mathbb{C})$ is an isomorphism of $G$-graded Lie algebras.

Let $P_n$ be the subgroup of $GL(n,\mathbb{C})$ of order $n^3$ generated by  $X_n$ and $Y_n$. Let  $K=\Int(\sl(n,\mathbb{C}))=PGL(n,\mathbb{C})$ and $\texttt{P}_n$ be the image of $P_n$ in  $\Int(\sl(n,\mathbb{C}))$, i.e., $\texttt{P}_n$ is generated by $Ad_X$ and $Ad_Y$. One knows that $$\wh{G}=\texttt{P}_n\cong \mathbb{Z}_n\times \mathbb{Z}_n.$$ It is proved that $\wh{G}$ is a maximal quasitorus of $\Aut(\sl(n,\mathbb{C}))$ in \cite{hpp}. So this is a good grading on $\sl(n,\mathbb{C})$. By \cite{a} one knows that $$W_K(\wh{G})\cong W_K(\Ga)=W(R)\cong SL(2,\mathbb{Z}_n).$$
\end{ex}

\begin{ex}\lb{b3}
 Assume $n=n_1\cdots n_k$ with each $n_i>1$ and \be \lb{c1} n_i|n_{i+1}, i=1,\cdots, k-1.\ee For $t=1,\cdots,k$, let $G_t=\mathbb{Z}_{n_t}\times \mathbb{Z}_{n_t}.$ Let $\bt_t:G_t\times G_t\rt \mathbb{C}^\times, \bt_t((i,j),(s,l))=\varepsilon_t^{js-il}$ be the alternating bicharacter on $G_t$ as in last example, where $\varepsilon_t$ is a primitive $n_t$-th root of unity. Let $(G,\bt)$ be the  orthogonal direct product symplectic abelian group of $(G_1,\bt_1),\cdots, (G_k,\bt_k)$, i.e.,\[G=G_1\times \cdots\times G_k\] and \be\lb{z}\bt((a_1,\cdots,a_k),(b_1,\cdots,b_k))=\bt_1(a_1,b_1)\cdots \bt_k(a_k,b_k)\ee where $a_i,b_i\in G_i$. Then $\bt$ is a nonsingular alternating bicharacter on $G$ and $R=G\setminus \{0\}$ is a finite root system in $(G,\bt)$. Let  $\xi\in Hom_\mathbb{Z}(G\otimes G, \mathbb{C}^\times)$ be defined by \[\xi(((i_1,j_1),\cdots,(i_k,j_k)),((s_1,t_1),\cdots,(s_k,t_k)))=\prod_{i=1}^k\varepsilon_i^{j_i s_i}.\] Then $\Psi(\xi)=\bt$.

 Next we will show that  $L(R)\cong \sl(n,\mathbb{C})$.

 Let $M=M(n,\mathbb{C})$. Then $M\cong M({n_1},\mathbb{C})\ot\cdots\ot M({n_k},\mathbb{C})$ has a $G$-grading as follows. Assume $ M({n_i},\mathbb{C})=\op_{(s,t)\in G_i} \mathbb{C} X_{n_i}^s Y_{n_i}^{-t}$ is the standard  $\varepsilon_i$-grading defined in last example. Then $M$ has the tensor product grading
 \be \lb{b5} M=\bigoplus_{(i_1,j_1,\cdots,i_k,j_k)\in G_1\times\cdots\times G_k} \mathbb{C}\cdot X_{n_1}^{i_1} Y_{n_1}^{-j_1}\otimes\cdots\otimes X_{n_k}^{i_k} Y_{n_k}^{-j_k}\ee
and $\sl(n,\mathbb{C})$ has the grading
 \be \lb{b6} \sl(n,\mathbb{C})=\bigoplus_{(i_1,j_1,\cdots,i_k,j_k)\in G_1\times\cdots\times G_k\setminus \{0\}} \mathbb{C}\cdot X_{n_1}^{i_1} Y_{n_1}^{-j_1}\otimes\cdots\otimes X_{n_k}^{i_k} Y_{n_k}^{-j_k}.\ee

 Define \[\varphi:L(R)\rt \sl(n,\mathbb{C}), u_{((i_1,j_1),\cdots,(i_k,j_k))}\mapsto X_{n_1}^{i_1} Y_{n_1}^{-j_1}\otimes\cdots\otimes X_{n_k}^{i_k} Y_{n_k}^{-j_k}.\] As in last example one can verify that $\varphi:L(R)\rt \sl(n,\mathbb{C})$ is an isomorphism of $G$-graded Lie algebras.

 We remark that for $R$ to be a finite root system in $(G,\bt)$ the condition (\ref{c1}) is unnecessary. But with this condition different finite root systems are not isomorphic.

As $L(R)$ is simple, the finite root system $R$ in this example is reduced and irreducible. We say that it is of type $I(n_1,\cdots,n_k)$, with each $n_i>1$ and $n_i|n_{i+1}\ \textrm{for}\ i=1,\cdots, k-1.$ For example, the finite root system in last example is of type $I(n)$.

Now we will describe the dual group $\wh{G}$ in $\Aut(L(R))$.

Let $D_n$ be the subgroup of $GL(n,\mathbb{C})$ consisting of all the diagonal matrices and $\texttt{D}_n\st PGL(n,\mathbb{C})$ be its quotient group in $PGL(n,\mathbb{C})$ under the canonical projection $GL(n,\mathbb{C})\rt PGL(n,\mathbb{C})$. Let $\texttt{P}_n\cong \mathbb{Z}_{n}\times \mathbb{Z}_{n}$ be the subgroup of $PGL(n,\mathbb{C})$ as defined in Example \ref{c2}. Assume $n=t n_1\cdots n_k$, $t\geq 1, n_i>1$, then \[ \texttt{D}_t
 \times \texttt{P}_{n_1}\times \cdots \times \texttt{P}_{n_k}\hookrightarrow PGL(n,\mathbb{C})\] by the adjoint action of $D_t\ot P_{n_1}\ot\cdots\ot P_{n_k}$ on $M(t,\mathbb{C})\ot M(n_1,\mathbb{C})\ot \cdots\ot M(n_k,\mathbb{C})\cong M(n,\mathbb{C})$. By
Theorem 3.2 of \cite{hpp}, any maximal quasitorus of $\PGL(n,\mathbb{C})$ is conjugate to one and only one of the $\texttt{D}_t
 \times \texttt{P}_{n_1}\times \cdots \times \texttt{P}_{n_k}$ with $n=t n_1\cdots n_k$, $t\geq 1, n_i>1$ and each $n_i$ dividing $n_{i+1}$. Thus any finite maximal
 quasitorus in  $\PGL(n,\mathbb{C})$ (one must have $t=1$) is
 conjugate to \be\lb{y1}Q(n_1,\cdots,n_k)=\texttt{P}_{n_1}\times \cdots \times \texttt{P}_{n_k}\ee  where $n=n_1\cdots n_k$, $n_i>1$ and each $n_i$ divides $n_{i+1}$. The dual group $\wh{G}$ of $G=\prod_{i=1}^k \mathbb{Z}_{n_i}\times \mathbb{Z}_{n_i}$ in this example is $Q(n_1,\cdots,n_k)$ in (\ref{y1}), which is a maximal quasitorus in  $K=\Int(\sl(n,\mathbb{C}))= PGL(n,\mathbb{C})$. The group $Q(n_1,\cdots,n_k)$ is also a maximal quasitorus in  $\Aut(\sl(n,\mathbb{C}))$ except the case $n_1=\cdots=n_k=2$.

One knows from \cite{h2} that $$W_K(\wh{G})\cong W_K(\Ga)=W(R)\cong Sp(G,\bt).$$
\end{ex}
\bigskip
\noindent \textbf{5.2.} Finite root systems $(G,\bt,R)$ such that $W(R)$ is a symmetric group

 Let $G_1=\mathbb{Z}_2^{n}$. For $i=1,\cdots, n$, let $e_i=(0,\cdots,0,1,0,\cdots,0)\in G_1$ with 1 in the $i$-th position. Let $R=\{e_i+e_j|1\leq i<j\leq n\}$ and $G\cong \mathbb{Z}_2^{n-1}$ be the subgroup of $G_1$ generated by $R$.

Let \[\bt:G_1\times G_1\rt \mathbb{C}^\times, \bt(a,b)=(-1)^{(\sum_{i=1}^{n} a_i)(\sum_{i=1}^{n} b_i)-\sum_{i=1}^{n} a_i b_i}.\] Let \[\xi:G_1\times G_1\rt \mathbb{C}^\times, \xi(a,b)=(-1)^{\sum_{1\leq j<i\leq n}(a_i b_j)}.\] Then $\Psi(\xi)=\bt$. We will show that $R$ is a finite root system in $(G,\bt)$.
Now we distinguish the cases $n$ is odd or even.

\begin{ex}\lb{v}
 Assume $n=2k+1$ is odd with $k\geq 1$.

Let $G_1=\mathbb{Z}_2^{2k+1}$. It is directly verified that $\Rad(\bt)=\{0,(1,\cdots,1)\}$.
Now $R=\{e_i+e_j|1\leq i<j\leq n=2k+1\}$ and $G\cong \mathbb{Z}_2^{2k}$ is the subgroup generated by $R$. It is easy to verify that $\bt|G$ is nonsingular, and that $R\subset G$ satisfies FRS0, FRS1 and FRS2, thus is a reduced finite root system in $(G,\bt)$. Then we will show that $L(R)\cong \so(2k+1,\mathbb{C})$.

For any $e_i+e_j\in R$ we will always assume $i<j$. Define \[\varphi:L(R)\rt \so(2k+1,\mathbb{C}), \varphi(u_{e_i+e_j})=2(E_{ij}-E_{ji}).\] Then for any $i,j,s,t\in \{1,\cdots,2k+1\}$ with $i<j,s<t$, by simple computations one has
\[\xi(e_i+e_j,e_s+e_t)-\xi(e_s+e_t,e_i+e_j)=0\ \textrm{if}\ \{i,j\}=\{s,t\}\ \text{or}\ \{i,j\}\cap\{s,t\}={\O}, \]
\[\xi(e_i+e_j,e_i+e_t)-\xi(e_i+e_t,e_i+e_j)=-2\ \textrm{if}\ j<t,\] and
\[\xi(e_i+e_j,e_j+e_t)-\xi(e_j+e_t,e_i+e_j)=2\ \textrm{if}\ j<t.\] Then
\[[u_{e_i+e_j},u_{e_s+e_t}]=\begin{cases}0 &\text{if $\{i,j\}=\{s,t\}$ or $\{i,j\}\cap\{s,t\}={\O};$ }     \\
-2 u_{e_j+e_t} &\text{if $i=s,j<t$;  }\\
2 u_{e_i+e_t} &\text{if $s=j<t.$  }
       \end{cases}\]
Also one has
\[[2(E_{ij}-E_{ji}),2(E_{st}-E_{ts})]=\begin{cases}0 &\text{if $\{i,j\}=\{s,t\}$ or $\{i,j\}\cap\{s,t\}={\O}$; }     \\
-2\cdot 2(E_{jt}-E_{tj}) &\text{if $i=s,j<t$; }\\
 2\cdot 2(E_{it}-E_{ti})&\text{if $s=j<t$.  }
       \end{cases}\]

 So $\varphi:L(R)\cong \so(2k+1,\mathbb{C})$ is an isomorphism of Lie algebras and the standard grading on $L(R)$ is in fact the  $\mathbb{Z}_2^{2k}$-grading $\Ga$ on $\so(2k+1,\mathbb{C})$: \[\Ga:L=\bigoplus_{1\leq i<j\leq 2k+1}L_{e_i+e_j} \]  with $L=\so(2k+1,\mathbb{C})$ and $L_{e_i+e_j}=\mathbb{C}(E_{ij}-E_{ji})$.

As $L(R)$ is simple, the finite root system $R$ in this case is reduced and irreducible. We say that it is of type $II(k)$ with $k\geq 1$. As $\sl(2,\mathbb{C})\cong \so(3,\mathbb{C})$, it is easy to see that \[II(1)\cong I(2).\]

In this paper let  \[SO(n,\mathbb{C})=\{ A\in GL(n,\mathbb{C})|A^t A=I_n\}.\]

Let $K=\Int(\so(2k+1,\mathbb{C}))=SO(2k+1,\mathbb{C})$. One knows that $\wh{G}\cong \mathbb{Z}_2^{2k}$ embeds as \be\lb{y2} SO(2k+1,\mathbb{C})\cap \{diag(\pm 1,\cdots,\pm 1\}),\ee thus is a maximal quasitorus in $K$. So this is a good grading on $\so(2k+1,\mathbb{C})$.

 By \cite{a} one has $$W_K(\wh{G})\cong W_K(\Ga)=W(R)\cong S_{2k+1}=S_n.$$

\end{ex}

\begin{ex}\lb{b7}
Assume $n=2k$ is even with $k\geq 3$.

(1) Let  $G_1=\mathbb{Z}_2^{2k}$. Then $\Rad(\bt)=\{0\}$. Now
$R=\{e_i+e_j|1\leq i<j\leq 2k\}$ and $G\cong \mathbb{Z}_2^{2k-1}$ is the subgroup generated by $R$. Then $\Rad(\bt|G)=\{0,(1,\cdots,1)\}$.
As in last example,
$R$ is a (nonreduced) finite root system in $(G,\bt)$.  And  \[\varphi:L(R)\rt \so(2k,\mathbb{C}), \varphi(u_{e_i+e_j})=2(E_{ij}-E_{ji})\] defines an isomorphism of $L(R)$ and $\so(2k,\mathbb{C})$, and the standard grading on $L(R)$ is in fact a $\mathbb{Z}_2^{2k-1}$-grading on $\so(2k,\mathbb{C})$.

As $\so(2k,\mathbb{C})$ is simple for $k\geq 3$ and  $\Rad(\bt|G)\neq 0$, the finite root system ${R}$ is nonreduced and irreducible. We say that it is of type $IV^{'}(k)$ with $k\geq 3$.

For $k\neq 4$, $\Aut(\so(2k,\mathbb{C}))\cong PO(2k,\mathbb{C})=O(2k,\mathbb{C})/\{\pm I_{2k}\}$. If $k=4$ then $PO(2k,\mathbb{C})$ is a subgroup of $\Aut(\so(2k,\mathbb{C}))$ of index 3.  One knows that $\wh{G}$ embeds as
\be \lb{y5} (O(2k,\mathbb{C})\cap \{diag(\pm 1,\cdots,\pm 1\})/\{\pm I_{2k}\}\ee in $PO(2k,\mathbb{C})$ which is
 a maximal quasitorus in $\Aut(\so(2k,\mathbb{C})$ (including the case $k=4$). So this is a good grading on $\so(2k,\mathbb{C})$.

Let $K=\Int(\so(2k,\mathbb{C}))$. Then by simple computation one has $$W_K(\wh{G})\cong W_K(\Ga)=W(R)\cong S_{2k}=S_n.$$

(2) (Continued) Let $H=\Rad(\bt|G)\cong \mathbb{Z}_2$ and $\overline{G}=G/H\cong \mathbb{Z}_2^{2k-2}$. Then as $\so(2k,\mathbb{C})$ is simple for $k\geq 3$, by Lemma \ref{l3}, $\overline{R}$ is a reduced and irreducible finite root system in $\overline{G}$. We say that it is of type $IV(k)$ with $k\geq 3$. As $\so(6,\mathbb{C})\cong \sl(4,\mathbb{C})$,it is not hard to check that \[IV(3)\cong I(2,2).\]

Let $T$ be the dual group of $\overline{G}$, which embeds as \be\lb{y3}(SO(2k,\mathbb{C})\cap \{diag(\pm 1,\cdots,\pm 1\})/\{\pm I_{2k}\}\ee in $\Aut(\so(2k,\mathbb{C}))=O(2k,\mathbb{C})/\{\pm I_{2k}\}$, which is a maximal quasitorus in $\Int(\so(2k,\mathbb{C}))=SO(2k,\mathbb{C})/\{\pm I_{2k}\}$. So this is a good inner grading on $\so(2k,\mathbb{C})$.

Let $K=\Int(\so(2k,\mathbb{C}))$. By \cite{a} one also has $$W_K(T)\cong W_K(\Ga)=W(R)\cong S_{2k}=S_n.$$

\end{ex}
\bigskip
\noindent \textbf{5.3.} Finite root systems constructed from nonsingular quadratic forms over $\mathbb{F}_2$

Let $q:\mathbb{F}_2^{2k+1}\rt \mathbb{F}_2$ be a  quadratic form and $\al_1:\mathbb{F}_2^{2k+1}\times \mathbb{F}_2^{2k+1}\rt \mathbb{F}_2$ be its polarization, i.e., $\al_1(a,b)=q(a+b)-q(a)-q(b)$. Assume that $q$ is nonsingular, i.e.,  $\al_1$ is an alternating bilinear form whose radical is 1-dimensional and $\Rad(q)=0$, where  $\Rad(q)=\{a\in \Rad(\al_1)|q(a)=0.\}$.  One knows that there is only one type of nonsingular quadratic form on $\mathbb{F}_2^{2k+1}$ up to isomorphism, and \be h:\mathbb{F}_2^{2k+1}\rt \mathbb{F}_2, h(a)=\sum_{i=1}^k a_{2i-1}a_{2i}+a_{2k+1}^2\ee is a such nonsingular quadratic form. One has \[\al_1(a,b)=h(a+b)-h(a)-h(b)=\sum_{i=1}^k (a_{2i}b_{2i-1}-a_{2i-1}b_{2i}) \] Identify the group $G=\mathbb{Z}_2^{2k+1}$ with the additive group of $\mathbb{F}_2^{2k+1}$. Let \[\al: G\times G\rt \mathbb{C}^\times, \al(a,b)=(-1)^{\al_1(a,b)}. \]

\begin{lem}\lb{x4}  Assume $k\geq 2$. Then the subset \[R=\{a\in G\setminus \Rad(\al)|h(a)=1\}\] is a finite root system in the symplectic abelian group $(G,\al)$.
\end{lem}
\bp
FRS1 is clear. Now we prove FRS2. Note that $\Rad(\al)=\Rad(\al_1)=\{ 0,e_{2k+1}\}$. Assume that $a,b\in R$ and $\al(a,b)\neq 1$, i.e., $\al_1(a,b)\neq 0$. Then $h(a+b)=h(a)+h(b)+\al_1(a,b)=1$. If $a+b=0$ then $a=b$ which contradicts to $\al(a,b)\neq 1$. If $a+b=e_{2k+1}$ then $h(a)\neq h(b)$ which contradicts to $a,b\in R$. Thus $a+b\in R$.

Finally we prove FRS0. It is easy to see that the subset \[B=\{e_{2i-1}+e_{2i},\  e_{2i}+e_{2k+1}|i=1,\cdots,k\}\cup \{\sum_{j=2k-3}^{2k+1}e_j\}  \] of $R$ generates $G$.
\ep
This finite root system $(\mathbb{Z}_2^{2k+1},\al, R)\ (k\geq 2)$ is nonreduced with $\Rad(\al)=\{0,e_{2k+1}\}$, and is said to be of type $I^{'}(2,\cdots,2)$, where the number of 2's is $k$.
\begin{lem}
   Assume $k\geq 2$. The corresponding reduced finite root system of $(\mathbb{Z}_2^{2k+1},\al, R)$ is of type $I(2,\cdots,2)$. One has $(\mathbb{Z}_2^{2k+1},\al, R)$ is nonreduced and irreducible, and $L(R)=L(\overline{R})=\sl(2^k,\mathbb{C}).$
\end{lem}

\bp Let $G=\mathbb{Z}_2^{2k+1}$.
As $\Rad(\al)=\{0,e_{2k+1}\}$, $\overline{G}=G/ \Rad(\al)= \mathbb{Z}_2^{2k}$, and the quotient group homomorphism is $p:\mathbb{Z}_2^{2k+1}\rt \mathbb{Z}_2^{2k}, (a_1,\cdots,a_{2k+1})\mapsto (a_1,\cdots,a_{2k})$. Thus  \[\overline{\al_1}(a,b)=\sum_{i=1}^k (a_{2i}b_{2i-1}-a_{2i-1}b_{2i}) \] and \[\overline{\al}:\overline{G}\times\overline{G}\rt \mathbb{C}^\times, \overline{\al}(a,b)=(-1)^{\overline{\al_1}(a,b)} \] is nonsingular.

Next we will show that $\overline{R}=\mathbb{Z}_2^{2k}\setminus\{0\}$. It is clear that $\overline{R}\subseteq \mathbb{Z}_2^{2k}\setminus\{0\}$. For any $a\in \mathbb{Z}_2^{2k} \setminus\{0\}$,
if $h(a,0)=1$ then $(a,0)\in R$ and $p(a,0)=a$; if $h(a,0)=0$ then $(a,1)\in R$ and $p(a,1)=a$. Thus the corresponding reduced finite root system of $(\mathbb{Z}_2^{2k+1},\al, R)$ is $(\mathbb{Z}_2^{2k},\overline{\al},\mathbb{Z}_2^{2k}\setminus\{0\})$, which is of type $I(2,\cdots,2)$. As $L(\overline{R})\cong \sl(2^k,\mathbb{C})$ is simple, by Lemma \ref{l3}, $R$ is irreducible and $L(R)=L(\overline{R})=\sl(2^k,\mathbb{C}).$
\ep

Recall that $\Aut(\sl(2^k,\mathbb{C}))=<PGL(2^k,\mathbb{C}),\tau>$  where $\tau:\sl(2^k,\mathbb{C})\rt \sl(2^k,\mathbb{C}), A\mapsto -A^t $ is an involutive outer automorphism of $\sl(2^k,\mathbb{C})$ for $k>1$. The dual group $\wh{G}$ is embedded as the maximal quasitorus \be \lb{y4}<Q(2,\cdots,2),\tau>\cong \mathbb{Z}_2^{2k+1}\ee in $\Aut(\sl(2^k,\mathbb{C}))$.

Similar as in Example \ref{b3}, one has $$W_K(\wh{G})\cong W_K(\Ga)=W(R)\cong Sp(2k,\mathbb{F}_2).$$
\bigskip

\begin{ex}\lb{bb9}
Assume $G=\mathbb{Z}_2^{2k}$ with $k\geq 3$. Let $R=\{a\in G|g(a)=1\}$.By Lemma \ref{b9} $R$ is a finite root system in $(G,\bt)$. One knows that
\[L(R)\cong \so(2^k,\mathbb{C})\] by Lemma \ref{a4} (1), and the corresponding grading is (\ref{z1}).

So the finite root system ${R}$ is reduced and irreducible.

On $\mathbb{C}^2$ there is a standard nonsingular orthogonal bilinear form $\varphi_0$ with matrix $I_2$. It is clear that $P_2\st O(2,\mathbb{C})$.  Equip $(\mathbb{C}^2)^{\ot k}\cong \mathbb{C}^{2k}$ with the nonsingular symmetric bilinear form $\varphi_0^{\ot k}$.  Then $P_2\ot\cdots \ot P_2\st SO(2^k,\mathbb{C})$ for $k\geq 2$. Then $Q(2,\cdots,2)=\texttt{P}_2\times\cdots\times \texttt{P}_2(\cong \mathbb{Z}_2^{2k})\st PSO(2^k,\mathbb{C})=\Int(\so(2^k,\mathbb{C}))$. It is directly verified that in this case the dual group $\wh{G}$ of $ G=\mathbb{Z}_2^{2k}$ embeds as the finite maximal quasitorus \be\lb{y7} Q(2,\cdots,2)=\texttt{P}_2\times\cdots\times \texttt{P}_2\cong \mathbb{Z}_2^{2k} \ee in $\Aut(\so(2^k,\mathbb{C}))$. We say that it is of type $V(k)$ with $k\geq 3$.

By \cite{a} one knows that $$W_K(\wh{G})\cong W_K(\Ga)=W(R)\cong O(G,g),$$ where $G=\mathbb{Z}_2^{2k}$ and
 $O(G,g)$ denotes the group of linear isomorphisms of $G$ preserving $g$, usually denoted by $O_{+}^{2k}(2)$.

 Note that the abelian groups in  $IV(4)$ and $V(3)$ are both $\mathbb{Z}_2^6$ and the respective $L(R)$ are both isomorphic to $\so(8,\mathbb{C})$, but
\begin{lem} \lb{aa5} \[IV(4)\ncong V(3).\]\end{lem}
\bp
Assume that $(G_1,\bt_1,R_1)$ is of the type $IV(4)$ and $(G_2,\bt_2,R_2)$ is of the type $V(3)$. Then $\wh{G_1}$ embeds in $\Aut(\so(8,\mathbb{C}))=PO(8,\mathbb{C})$ as in (\ref{y3}), which is not maximal abelian in $PO(8,\mathbb{C})$ as its centralizer in $PO(8,\mathbb{C})$ is as in (\ref{y5}). Since $\wh{G_2}$ is a maximal abelian subgroup in $\Aut(\so(8,\mathbb{C}))$, $\wh{G_1}$ and $\wh{G_2}$ cannot be conjugate in $\Aut(\so(8,\mathbb{C}))$ and $IV(4)\ncong V(3)$.
\ep
\end{ex}

\begin{ex}\lb{b8}
Assume  $G=\mathbb{Z}_2^{2k}$ with $k\geq 1$. Let $R=\{a\in G|f(a)=1\}$. By Lemma \ref{bb8} $R$ is a finite root system in $(G,\bt)$.  One knows that
\[L(R)\cong\sp(2^k,\mathbb{C})\] by Lemma \ref{a4} (2), and the corresponding grading is (\ref{z2}).

So the finite root system ${R}$ is reduced and irreducible. We say that it is of type $III(k)$ with $k\geq 1$. As $\sp(2,\mathbb{C})\cong \sl(2,\mathbb{C})$ and $\sp(4,\mathbb{C})\cong \so(5,\mathbb{C})$, one can verify that \[III(1)\cong I(2)\  and\  III(2)\cong II(2).\]

On $\mathbb{C}^2$ there is a standard nonsingular skew-symmetric bilinear form $\varphi_1$ with matrix $\left(
             \begin{array}{cc}
               0 & 1\\
               -1 & 0 \\
             \end{array}
           \right)$. Equip $\mathbb{C}^2\ot\cdots\ot \mathbb{C}^2\cong \mathbb{C}^{2k}$ with the nonsingular skew-symmetric bilinear form $\varphi_1\ot\varphi_0\ot\cdots\ot\varphi_0$. As $P_2\st \mathbb{C}^\times Sp(2,\mathbb{C})$ and $P_2\st O(2,\mathbb{C})$, $P_2\ot\cdots \ot P_2$  acts on $\mathbb{C}^2\ot\cdots\ot \mathbb{C}^2\cong \mathbb{C}^{2k}$ and $P_2\ot\cdots \ot P_2\st \mathbb{C}^\times Sp(2^k,\mathbb{C})$. Then $Q(2,\cdots,2)=\texttt{P}_2\times\cdots\times \texttt{P}_2(\cong \mathbb{Z}_2^{2k})\st PSp(2^k,\mathbb{C})=\Aut(\sp(2^k,\mathbb{C}))$. In this case the dual group $\wh{G}$ of $ G=\mathbb{Z}_2^{2k}$ embeds as the finite maximal quasitorus \be\lb{y8} Q(2,\cdots,2)=\texttt{P}_2\times\cdots\times \texttt{P}_2\cong \mathbb{Z}_2^{2k} \ee in $\Aut(\sp(2^k,\mathbb{C}))$.

By \cite{a} one knows that $$W_K(\wh{G})\cong W_K(\Ga)=W(R)\cong O(G,f),$$ where $G=\mathbb{Z}_2^{2k}$ and
 $O(G,f)$ denotes the group of linear isomorphisms of $G$ preserving $f$, usually denoted by $O_{-}^{2k}(2)$.
\end{ex}

\begin{prop}\lb{aa6}
Assume that  $L$ is a complex simple Lie algebra, $T$ is a good finite maximal quasitorus in $\Aut(L)$ or in $\Int(L)$. If  $(L,T)$  is as in above examples of this section, then there corresponds to $(L,T)$ a unique  finite root system $(G,\bt,R)$ up to isomorphism such that $L=L(R)$ and the standard grading on $L(R)$ is isomorphic to the grading on $L$ induced by $T$.
\end{prop}
\bp
 In the above examples, we constructed finite root systems $(G,\bt,R)$ and identified $L=L(R)$ with the corresponding standard gradings on it, then gave the embedding of $T=\wh{G}$ in $\Aut(L)$ or $\Int(L)$ and showed that $T$ is a finite maximal quasitorus in $\Aut(L)$ or $\Int(L)$. So $(G,\bt,R)$ is good, and is clearly a finite root system corresponding to $(L,T)$. We only need to prove that it is the unique  finite root system corresponding to $(L,T)$ up to isomorphism.

If $(L,T)$ is as in Example \ref{b3} (Example \ref{c2} is a special case of Example \ref{b3}), then $L=\sl(n,\mathbb{C})$, $T$ is as in \ref{y1} and the grading is as in \ref{b6}. Thus if  $(G,\bt,R)$ is a finite root system corresponding to $(L,T)$, then $G=\wh{T}$ and $R=G\setminus\{0\}$. Since $R\subseteq G\setminus \Rad(\bt)$, one has $\Rad(\bt)=\{0\}$ and thus $\bt$ is nonsingular. By  Theorem 1.8 of \cite{ka3}, if there are two nonsingular alternating bicharacters $\bt_i $ on $G$ for $i=1,2$, then $(G,\bt_1)$ and $(G,\bt_2)$ are isometric. So the finite root system $(G,\bt,R)$ corresponds to $(L,T)$ is  unique up to isomorphism.

 If $(L,T)$ is as in other examples, then $G=\wh{T}$ is always an elementary abelian 2-groups. So the finite root system corresponding to $(L,T)$ is unique up to isomorphism by Lemma \ref{x7}.

\ep

\section{Classification of good finite maximal quasitorus in $\Aut(L)$ with $L$ a complex simple Lie algebra and main results}
 \setcounter{equation}{0}\setcounter{theorem}{0}
Now we state our main result.

\begin{theorem}\lb{b1}
Assume that $L$ is a complex simple Lie algebra, $T$ is a good finite maximal quasitorus in $\Aut(L)$ or in $\Int(L)$.

(1) There corresponds to $(L,T)$ a unique irreducible finite root system $(G,\bt,R)$ up to isomorphism.

(2) The pair $(L,T)$ and corresponding finite root system $R$, where $T$ is up to conjugation in $\Aut(L)$, is one and only one of the following.

1. $L=\sl(n,\mathbb{C}), n\geq 2$.

(a) Assume $n=n_1\cdots n_k$ with $n_i|n_{i+1}$ for $i=1,\cdots,k-1$ and $n_i>1$. Then $T\cong \mathbb{Z}_{n_1}^2\times\cdots\times \mathbb{Z}_{n_k}^2$ which embeds in $\Int(L)$ as in (\ref{y1}) is always a good finite maximal quasitorus in  $\Int(L)$, which is also a good  finite maximal quasitorus in $\Aut(L)$ except $k\geq 2$ and $n_1=\cdots=n_k=2$. The corresponding finite root system $R$ is of type $I(n_1,\cdots,n_k)$.

(b) Assume $n=2^k, k\geq 2$. Then $T\cong \mathbb{Z}_2^{2k+1}$ which embeds in $\Aut(L)$ (and $T\nsubseteq \Int(L)$) as in (\ref{y4}) is a good  finite maximal quasitorus in $\Aut(L)$. The corresponding finite root system $R$ is of type $I^{'}(2,\cdots,2)$.

2.  $L=\so(2k+1,\mathbb{C}),k\geq 2$. Then $T\cong \mathbb{Z}_2^{2k}$ which embeds in $\Aut(L)$ as in (\ref{y2}) is a good  finite maximal quasitorus in $\Aut(L)$. The finite root system $R$ is of type $II(k)$. In this case $\Aut(L)=\Int(L)$ and $T$ is also a good finite maximal quasitorus in $\Int(L)$.

3. $L=\sp(2^k,\mathbb{C}),k\geq 3$, Then $T\cong \mathbb{Z}_2^{2k}$ which embeds in $\Aut(L)$ as in (\ref{y8}) is a good  finite maximal quasitorus in $\Aut(L)$. The finite root system $R$ is of type $III(k)$. In this case $\Aut(L)=\Int(L)$ and $T$ is also a good finite maximal quasitorus in $\Int(L)$.

4. $L=\so(2k,\mathbb{C}),k\geq 4$. Then $T\cong \mathbb{Z}_2^{2k-1}$ which embeds in $\Aut(L)$ (and $T\nsubseteq \Int(L)$) as in (\ref{y5}) is a good  finite maximal quasitorus in $\Aut(L)$. The finite root system $R$ is of type $IV^{'}(k)$. In this case $T\cap Int(L)\cong \mathbb{Z}_2^{2k-2}$is a good  finite maximal quasitorus in $\Int(L)$, and the finite root system $R$ is of type $IV(k)$.

5. $L=\so(2^k,\mathbb{C}),k\geq 3$. Then $L\cong \mathbb{Z}_2^{2k}$ which embeds in $\Int(L)$ as in (\ref{y7}) is a good  finite maximal quasitorus in $\Aut(L)$ as well as in $\Int(L)$, and the finite root system $R$ is of type $V(k)$.

In particular there is no good finite maximal quasitorus in $\Aut(L)$ for exceptional simple Lie algebras $L$.
\end{theorem}
We will prove it later.
\begin{rem}\lb{aa1}
If $T$ is a good finite maximal quasitorus in $\Aut(L)$ for $L$ a simple Lie algebra, then either $T\subseteq \Int(L)$ or $T\nsubseteq \Int(L)$, in the latter case $T\cap\Int(L)$ is a  good finite maximal quasitorus in $\Int(L)$. The corresponding irreducible finite root system of $(L,T)$ is reduced if $T$ is a good finite maximal quasitorus in $\Int(L)$, and is nonreduced if $T$ is a good finite maximal quasitorus in $\Aut(L)$ but $T\nsubseteq\Int(L)$.
\end{rem}
\begin{rem}
Let $K=\Int(L)$ and $T=\wh{G}$. One knows from \cite{a} and \cite{h2} that in all the cases in the list, the Weyl group $W_K(T)$ of $T$ with respect to $K$ is isomorphic to the Weyl group $W(R)$, i.e.,  $W_K(T)$  is generated by root transvections. So Vogan's conjecture holds in these cases. The Weyl group in Case (a) is $Sp(G,\bt)$ by \cite{h2}. And the Weyl group in Case (b), (c), (d) and (e) are respectively $S_{2k+1}$, $O(G,f)$, $S_{2k}$ and $O(G,g)$ by \cite{a}, where $O(G,f)$ and $O(G,g)$ are the respective group of linear isomorphisms of $G$ preserving $f$ and $g$.  $O(G,f)$ and $O(G,g)$ are usually denoted by $O_{-}^{2k}(2)$ and  $O_{+}^{2k}(2)$ respectively. One also notes in \cite{a} that in Case (b), (c), (d) and (e), $G$ is a subgroup of $SL(n,\mathbb{F}_2)$ acting irreducible on $\mathbb{F}_2^n$.
\end{rem}
\begin{coro}
Assume that $(G,\bt,R)$ is a finite root system such that $L(R)$ is a simple Lie algebra and the standard grading on  $L(R)$ is a good inner grading, then $(G,\bt,R)$ is reduced and irreducible, and the type of $(G,\bt,R)$ is one and only one of the following:

(1) $I(n_1,\cdots,n_k)$ with $k\geq 1$, each $n_i>1$ and $n_i|n_{i+1}$ for $i=1,\cdots,k-1$;

(2) $II(k)$ with $k\geq 2$;

(3) $III(k)$ with $k\geq 3$;

(4) $IV(k)$ with $k\geq 4$;

(5) $V(k)$ with $k\geq 3$.
\end{coro}

\bp
By Theorem \ref{b1} and Remark \ref{aa1}, one only need to verify that there is no redundancy in the list. Assume that two finite root systems $(G_i,\bt_i,R_i)$ ($i=1,2$) in the list are isomorphic. Then $L(R_1)\cong L(R_2)$, which is possible only if they are of types $IV(m)$ and $V(k)$ respectively, where $m\geq 4, k\geq 3$ and \be\lb{aa3}2m=2^k.\ee Assume that $(G_1,\bt_1,R_1)$ is of type $IV(m)$ and $(G_2,\bt_2,R_2)$ is of type $V(k)$. Then $\wh{G_1}\cong \mathbb{Z}_2^{2m-2}$ and $\wh{G_2}\cong \mathbb{Z}_2^{2k}$ are isomorphic. So\be\lb{aa4}2m-2=2k.\ee Solving (\ref{aa3}) and (\ref{aa4}) one has $m=4,k=3$. But $IV(4)\ncong V(3)$ by Lemma \ref{aa5}.

\ep

Now we will prove Theorem \ref{b1}. Recall that for a complex simple Lie algebra $L$, a finite maximal quasitorus $T$ of $\Aut(L)$ or $\Int(L)$ is said to be good if $\dim L_{\alpha}\leq 1$ for every character $\alpha\in \wh{T}$. The main part of the proof is to classify good finite maximal quasitorus $T$ of $\Aut(L)$ (or $\Int(L)$). If the classification is done, then for each $(L,T)$ the finite root system $(G,\bt,R)$ is constructed in Section 5. Theorem \ref{b1} (1) follows from Proposition \ref{aa6}.

Now we classify good finite maximal quasitorus $T$ of $\Aut(L)$ (or $\Int(L)$).
Since the finite  maximal quasitorus $T$ of $\Aut(L)$ will stabilize some compact
real form $L_0$ of $L$, $T\subset \Aut(L_0)$. Thus it is equivalent to the
classification of finite maximal quasitorus  of $\Aut(L_0)$ (or $\Int(L_0)$). To do
this we use the classification of maximal abelian subgroups of compact simple Lie
groups in \cite{y}. We will also show that there are no good gradings on exceptional simple
Lie algebras.

When $L_0=\mathfrak{su}(n)$ ($n\geq 2$) and $T\subset\Int(L_0)$ is a finite maximal
abelian subgroup, there is no much to say, as every finite maximal abelian subgroup is of the form (\ref{y1}) and satisfies
the condition. We remark that, in this case $T$ is also a maximal abelian
subgroup of $\Aut(L_0)$ except when it is an elementary abelian 2-subgroup.
In the latter case $T$ commutes with an outer involution $\tau$ in the conjugacy
class of complex conjugation and $\langle\tau,T\rangle$ is a maximal
abelian subgroup of $\Aut(L_0)$.

When $L_0=\mathfrak{so}(n)$ ($n\geq 5$) and $T\subset\PO(n)$ is a finite maximal
abelian subgroup(this means only the case of $L_0=\mathfrak{so}(8)$ and
$T/T\cap\Int(L_0)\cong C_3$ needs a further investigation), let
$\pi:\O(n)\rightarrow\PO(n)$ be the natural projection.
By \cite[Subsection 3.2]{y}, there is an antisymmetric and bimultiplicative
function $m:T\times T\rightarrow\{\pm{1}\}$ and $\ker m$ is a diagonalizable
subgroup (since $T$ is assumed to a finite maximal abelian subgroup,
$\ker m=B_{T}$ by \cite[Lemma 3.3]{y}, where $B_{T}$ is a subgroup defined in
\cite{y}). By \cite[Proposition 3.3]{y}, there exists integers $s_{0}\geq 1$
and $k\geq 0$ such that $n=s_{0}\cdot 2^{k}$, and the centralizer of $\pi^{-1}(T)$ in $\O(n)$ is
$$C_{\O(n)}(\pi^{-1}(T))=\underbrace{\O(2^{k})\times\cdots\times\O(2^{k})}_{s_0},$$
and $\pi^{-1}(\ker m)=Z(C_{\O(n)}(\pi^{-1}(T))$, the center of $C_{\O(n)}(\pi^{-1}(T))$. We show that $k=0$, or $s_0=1$,
or $(s_0,k)=(3,1)$. Moreover, there is a unique conjugacy class of
$T$ while $(s_0,k)=(3,1)$.
Let $E_{i}=\pi(\diag\{(i-1)\cdot I_{2^{k}},
I_{2^{k}},(s_{0}-i)\cdot I_{2^{k}}\})$. For any character $\alpha\in
\Hom(T,S^1)$, if $L_{\alpha}\neq 0$ and $\alpha|_{\ker m}
\neq 1$, then there exists $i\neq j$, $\alpha(E_{k})=-1$ if $k\neq i,j$
and $\alpha(E_{k})=1$ if $k= i$ or $j$. Given $i\neq j$ and any
character $\alpha\in\Hom(T,S^1)$ with $\alpha(E_{k})=-1$ if $k\neq i,j$
and $\alpha(E_{k})=1$ if $k= i$ or $j$, by calculation one shows
that $\dim L_{\alpha}=1$. Given a character $\alpha\in\Hom(T,S^1)$
with $\alpha|_{\ker m}=1$, the root space $\dim L_{\alpha}$ is
contained in the complexified Lie algebra of $$C_{\O(n)}(\pi^{-1}(T))=
\underbrace{\O(2^{k})\times\cdots\times\O(2^{k})}_{s_0}.$$ In \cite{y}, we
defined functions $\mu_1,\dots,\mu_{s_0}: T/\ker m\rightarrow\{\pm{1}\}$
by $$(A_1^{2},\dots,A_{s_0}^{2})=(\mu_{1}(x)I_{2^{k}},\dots,
\mu_{s_0}(x)I_{2^{k}})$$ for any $x\in T$ and $\diag\{A_1,\dots,A_{s_0}\}
\in\pi^{-1}(x)$. Let $L_{i}$ be the complexified Lie algebra
of the $i$-th factor of $C_{\O(n)}(\pi^{-1}(T))$. By calculation one
shows that $\dim(L_{\alpha}\cap L_{i})\leq 1$
and equality holds if and only if $\mu_{i}(\ker(m|_{\ker\alpha}))=-1$.
From these, we conclude that in the case of $s_0=1$ or $k=0$,
$\dim L_{\alpha}\leq 1$ for any character $\alpha$; in the case
of $s_0\geq 2$ or $k\geq 1$, $\dim L_{\alpha}\leq 1$ for any
character $\alpha$ if and only if for any $i\neq j$, there exists no
$x\in T$ with $\mu_{i}(x)=\mu_{j}(x)=-1$. Since $\{x\in T|\mu_{i}(x)=
\mu_{j}(x)\}$ is an index 2 subgroup, there exits no $x\in T$ with
$\mu_{i}(x)=\mu_{j}(x)=-1$ implies $k=1$ and $\mu_{i}\neq\mu_{j}$.
Moreover, we have $s_0=3$ since $n=s_{0}\cdot 2^{k}\geq 5$. Finally,
as $\mu_{1}$, $\mu_{2}$, $\mu_{3}$ are non-equal to each other, by
\cite[Proposition 3.4]{y}, we get a unique conjugacy class in the
case of $(s_0,k)=(3,1)$. As $\mathfrak{so}(6)\cong\mathfrak{su}(4)$,
actually one shows that the above subgroup $T$ of
$\Aut(\mathfrak{so}(6))$ corresponds to a
subgroup of $\Int(\mathfrak{su}(4))$ isomorphic to $(\mathbb{Z}/4\mathbb{Z})^{2}$ in the
case $(s_0,k)=(3,1)$. If $s_0=1$ then $n=2^k$ and the finite maximal quasitorus $T\cong \mathbb{Z}_2^{2k}$ is of the form (\ref{y7});
if $k=0$ then the finite maximal quasitorus $T\cong \mathbb{Z}_2^{n-1}$, which is of the form (\ref{y2}) if $n$ is odd and is of the form (\ref{y5}) if $n$ is even.

When $L_0=\mathfrak{sp}(n)$ ($n\geq 2$) and $T\subset\PSp(n)$ is a finite
maximal abelian subgroup,
similarly there are integers $s_{0}\geq 1$ and $k\geq 1$ with
$n=s_{0}\cdot 2^{k-1}$. An argument similar as the above for projective
orthogonal groups shows $s_0=1$ if $\dim L_{\alpha}\leq 1$ for
each $\alpha\in\Hom(T,S^{1})$. The case $k=0$ can not happen since
we have $k\geq 1$ automatically, and there is no exception like $(s_0,k)=(3,1)$
in this case. In the case $s_0=1$ the finite maximal quasitorus $T\cong \mathbb{Z}_2^{2k}$ is of the form (\ref{y8}).

When $L_0=\mathfrak{su}(n)$ ($n\geq 3$),
$T\subseteq\Aut(\mathfrak{su}(n))$ is a finite
maximal abelian subgroup, and $T\nsubseteq\Int(\mathfrak{su}(n))$,
similarly there are integers $s_{0}\geq 1$ and $k\geq 0$ with
$n=s_{0}\cdot 2^{k}$ (\cite[Proposition 3.10]{y}). Similarly as for
projective orthogonal groups, one shows that $s_0=1$ or $k=0$ if
$\dim L_{\alpha}\leq 1$ for each $\alpha\in\Hom(T,S^{1})$.
When $s_0=1$, $T$ is conjugate to a subgroup of the form
$\langle\tau_0, T'\rangle$, where $\tau_0$ is complex conjugation and
$T'$ is an elementary abelian $2$-subgroup of $\PO(2^{k})$ with rank
$2k$ and with a nondegenerate skew-symmetric and bimultiplicative
function $m: T'\times T'\rightarrow\{\pm{1}\}$. This subgroup
satisfies $\dim L_{\alpha}\leq 1$ for each weight $\alpha\in
\Hom(T,S^{1})$. Note that in this case the root spaces of the group grading induced by
$T$ and $T'$ coincides (that means if $\frg_{\alpha},\frg_{\beta}\neq 0$,
then $\alpha=\beta$ if and only if $\alpha|_{T'}=\beta|_{T'}$). The group $T'$
is a maximal elementary abelian $2$-subgroup of $\Int(\mathfrak{su}(n))$. The finite maximal quasitorus $T$ in the case $s_0=1$ is of the form (\ref{y4}).
When $k=0$, $T$ is conjugate to a subgroup of the form
$\langle\tau_0, T'\rangle$, where $\tau_0$ is complex conjugation and
$T'$ is an elementary abelian $2$-subgroup consisting of all diagonal
elements in $\PO(n)$. For the character with $\alpha|_{T'}=1$ and
$\alpha(\tau_0)=-1$, $\dim  L_{\alpha}=n-1>1$. Hence, this
subgroup is not good.

Now we turn to exceptional groups. If a finite abelian subgroup $T$ is
good, one must have $|T|\geq\dim L+1$. When
$L_0=\mathfrak{so}(8)$ and $T/T\cap\Int(L_0)
\cong C_3$, by \cite[Proposition 5.1]{y}, there are two conjugacy classes
of finite maximal abelian subgroups, with order $27$ and $24$ respectively.
Both orders are smaller than $\dim\frg=28$ and hence they are not good.

When $L_0$ is of type $\bf G_2$ (or $\bf F_4$), there is a
unique conjugacy classes (or two conjugacy classes) of finite maximal
abelian subgroups with order $8$ (or $27$, $32$), they are not good
as the order is less than $\dim\frg$.

When $L_0$ is of type $\bf E_6$, any finite maximal abelian
subgroup $T$ with $|T|>\dim L$ is conjugate to one of $F_3$,
$F_{12}$ (cf. \cite[Table 5 and Table 6]{y}). The group $F_3\cong(C_3)^{4}$
is not good since $\dim L_{\alpha}=3$ if
$\alpha|_{\langle\{x\in F_{3}| x\sim\theta_2\}\rangle}=1$. The
$F_{12}$ is not good since $\dim L_{\alpha}=6$ if
$\alpha|_{F_{12}\cap\Int(L_0)}=1$.

When $L_0$ is of type $\bf E_7$, any finite maximal abelian
subgroup $T$ with $|T|>\dim L$ is conjugate to $F_{7}$ (see
\cite[Table 7]{y}). The group $F_7\cong (C_2)^{8}$ is not good since $\dim
L_{\alpha}=7$ if $\alpha|_{F_{7}\cap\Int(L_0)}
=1$ (see \cite[Paragraph after Proposition 8.2]{y}).

When $L_0$ is of type $\bf E_8$, any finite maximal abelian
subgroup $T$ with $|T|>\dim L$ is conjugate to one of $F_4$,
$F_5$, $F_{7}$, $F_8$, $F_9$ (see \cite[Table 8]{y}). The group
$\Aut(\mathfrak{e}_8)$ has a Klein four-subgroup $\Gamma_1$ with
(see \cite[Table 6]{hy}) \[G^{\Gamma_1}\cong(\E_6\times\U(1)\times
\U(1))/\langle(c,\omega,1)\rangle\rtimes\langle z\rangle,\] where
$(\fre_6\oplus i\bbR\oplus i\bbR)^{z}=\frf_4\oplus 0\oplus 0$.
By \cite[Proposition 11.2]{y}, each of $F_4$, $F_5$, $F_{7}$, $F_8$ is
conjugate to a subgroup of the form $T=T'\times\Gamma_1$, where
$T'\subset\E_6\rtimes\langle z\rangle$. Since $T'$ is not good, $T$ is
not as well. The group $F_{9}$ is not good since $\dim L_{\alpha}=8$
if $\alpha|_{F_{9}\cap\Int(L_0)}=1$.

\end{document}